\documentclass{article}
\title{On Biproducts and Extensions}
\author{Yevgenia Kashina, Yorck Sommerh\"auser}
\date{}

\usepackage{amsmath}
\usepackage{amsthm}
\usepackage{amssymb}
\usepackage{mathtools}
\usepackage{needspace}

\makeatletter
\renewcommand{\subsection}{\@startsection{subsection}{2}{0em}%
{\baselineskip}{-0em}{\bfseries\normalsize}}
\makeatother

\newcounter{num}
\newenvironment{pflist}{\begin{list}{(\arabic{num})}{\usecounter{num} \leftmargin0cm \itemindent5pt}}{\end{list}}

\newenvironment{parlist}{\begin{list}{(\arabic{num})}{\usecounter{num} \leftmargin0pt \itemindent5pt \topsep0pt \labelwidth0pt}}{\end{list}}

\newenvironment{parlist1}{\begin{list}{(\arabic{num})}{\usecounter{num} \leftmargin0pt \itemindent5pt \itemsep7pt \topsep0pt \labelwidth0pt}}{\end{list}}

\newenvironment{thmlist}{\begin{list}{\arabic{num}.}{\usecounter{num} \leftmargin0pt \itemindent5pt \topsep-5pt \partopsep0pt \labelwidth0pt}}{\end{list}}

\newcounter{num1}
\newenvironment{proplist}{\begin{list}{\arabic{num1}.}{\usecounter{num1} \listparindent0pt \topsep-2pt  \leftmargin25pt \itemindent0pt \itemsep5pt \partopsep0pt \labelwidth10pt \labelsep0.5em }}{\end{list}}

\numberwithin{equation}{section}

\makeatletter
\def\enumfix{%
\if@inlabel
 \noindent \par\nobreak\vskip-\topsep\hrule\@height\z@
\fi}
\let\oldenumerate\enumerate
\def\enumerate{\enumfix\oldenumerate}
\let\oldthmlist\thmlist
\def\thmlist{\enumfix\oldthmlist}
\makeatother

\theoremstyle{definition}
\newtheorem*{defn}{Definition}
\newtheorem*{prop}{Proposition}
\newtheorem*{lem}{Lemma}
\newtheorem*{cor}{Corollary}
\newtheorem*{pf}{Proof}
\newtheorem*{thm}{Theorem}

\newtheoremstyle{break}{}{}{}{}{\bfseries}{.}{\newline}{}

\theoremstyle{break}

\newcommand{\dm}{\operatorname{dim}}

\newcommand{\Ind}{\operatorname{Ind}}
\newcommand{\Span}{\operatorname{Span}}

\newcommand{\Hom}{\operatorname{Hom}}
\newcommand{\id}{\operatorname{id}}
\newcommand{\co}{\scriptstyle \operatorname{co}}
\newcommand{\ot}{\mathbin{\otimes}}
\newcommand{\deq}{\vcentcolon =}
\newcommand{\sbb}{S_{B}}
\newcommand{\ua}{1_A}
\newcommand{\ub}{1_B}
\newcommand{\K}{1}

\newcommand{\xu}{x'}
\newcommand{\yu}{y'}
\newcommand{\uu}{u'}
\newcommand{\vu}{v'}
\newcommand{\ru}{r'}
\newcommand{\su}{s'}
\newcommand{\au}{a'}
\newcommand{\omu}{\omega'}
\newcommand{\etu}{\eta'}
\newcommand{\cu}{c'}
\newcommand{\du}{d'}
\newcommand{\hu}{h'}
\newcommand{\Au}{A'}
\newcommand{\Bu}{B'}

\newcommand{\1}{{(1)}}
\newcommand{\2}{{(2)}}
\newcommand{\3}{{(3)}}

\newcommand{\ra}{\rightarrow}

\newcommand{\da}{\Delta_{A}}
\newcommand{\db}{\Delta_{B}}

\newcommand{\ea}{{\varepsilon_{A}}}
\newcommand{\eb}{{\varepsilon_{B}}}
\newcommand{\eh}{{\varepsilon_{H}}}
\newcommand{\el}{{\varepsilon}}

\newcommand{\sa}{S_{A}}
\newcommand{\sd}{S_{C}}
\newcommand{\sh}{S_{H}}
\newcommand{\A}{1_A}
\newcommand{\B}{1_B}
\newcommand{\Z}{{\mathbb Z}}
\newcommand{\HH}{1_H}
\newcommand{\dhh}{\Delta_{H}}

\begin{document}

\maketitle

\begin{abstract}
\noindent
We describe in which ways the Radford biproducts of certain eight-dimen\-sional Yetter-Drinfel'd Hopf algebras over the elementary abelian group of order~$4$ can be written as extensions of Hopf algebras.
\end{abstract}


\section*{Introduction} \label{Sec:Introd}
In their article~\cite{KaSo2}, the authors described two semisimple Yetter-Drinfel'd Hopf algebras of dimension~$8$ over the group ring of the elementary abelian group of order~$4$. The purpose of the construction of these Yetter-Drinfel'd Hopf algebras was to show that the core of a group-like element is not always completely trivial. It was not even mentioned in this article that, via the Radford biproduct construction (cf.~\cite{RadfProj}), these Yetter-Drinfel'd Hopf algebras give rise to ordinary Hopf algebras of dimension~$32$. The purpose of the present article is to discuss how these Hopf algebras fit into the general theory of semisimple Hopf algebras of prime power dimension, as found for example in~\cite{Ka}, \cite{Ka2}, \cite{Ka3}, \cite{Ma1}, and~\cite{Ma2}. As we will see, these ordinary Hopf algebras also behave differently from other known examples.

Every semisimple Hopf algebra of prime power dimension can in principle be constructed as an iterated crossed product (cf.~\cite{MW}, Theorem~3.5). However, for the dimensions where the description is more concrete, or even so concrete that all semisimple Hopf algebras of that given dimension can be classified, like in the articles just cited, the Hopf algebras often contain a commutative Hopf subalgebra of prime index, a situation sometimes also encountered for semisimple Hopf algebras that are not of prime power dimension, like in~\cite{Ma3} and~\cite{Na1}. The examples constructed here do not contain such a Hopf subalgebra, as we show in Paragraph~\ref{Proof2by16}. More generally, we show that they cannot be constructed as central or cocentral abelian extensions.

Let us state the main results of the article more precisely, while we simultaneously explain its structure. Section~\ref{Sec:Prelim} contains preliminaries on Yetter-Drinfel'd Hopf algebras, Radford biproducts, and extensions of Hopf algebras. In Section~\ref{Sec:BiprodFirst}, we consider the eight-dimensional Yetter-Drinfel'd Hopf algebra~$A$ over the elementary abelian group of order~$4$ that appears in~\cite{KaSo2}, Section~2. We describe the arising Radford biproduct~$B$, which has dimension~$32$, give a presentation of~$B$ in terms of generators and relations, and compute its center. In Section~\ref{Sec:Rep}, we show that the groups of group-like elements of both~$B$ and~$B^*$ are elementary abelian of order~$8$, which means in particular that they have eight one-dimensional representations. In Theorem~\ref{Wedderb}, we show that, up to isomorphism, both~$B$ and~$B^*$ have in addition two irreducible representations of dimension~$2$ and one irreducible representation of dimension~$4$.

In Section~\ref{Sec:2by16}, we show that~$B$ contains a unique sixteen-dimensional Hopf subalgebra~$N$. It is normal and isomorphic to the Hopf algebra~$H_{d:1,1}$ from the classification appearing in~\cite{Ka}, Table~1. Our Hopf algebra~$B$ therefore fits into exactly one extension of the type
\[N \hookrightarrow B \twoheadrightarrow Z\]
with~$\dm N = 16$ and~$\dm Z = 2$. This extension is not abelian.

In Section~\ref{Sec:Dim2}, we show that the Hopf subalgebra~$N$ also arises as the unique sixteen-dimensional quotient of~$B$. By construction, $N$ is also a Radford biproduct, and the quotient map \mbox{$\pi_N \colon B \to N$} is induced from a map between the corresponding Yetter-Drinfel'd Hopf algebras. Although $\pi_N$ does not restrict to the identity on~$N$, it is conormal, and therefore~$B$ fits into exactly one extension of the type
\[U \hookrightarrow B \overset{\pi_N}{\twoheadrightarrow} N\]
with~$\dm U = 2$ and~$\dm N = 16$. This extension is not abelian. We use these facts to describe the structure of the Grothendieck ring of~$B$ in Paragraph~\ref{GrothB}.

In Section~\ref{Sec:4by8}, we begin by determining the eight-dimensional Hopf subalgebras of~$B$. There are three, denoted by~$M_1$, $M_2$, and~$M_3$, all of which are contained in~$N$. The algebra~$M_1$ is just the span of the group-like elements, while the other two are dual to the group ring of the dihedral group. We show in Paragraph~\ref{normal8} that, of the three Hopf subalgebras, only~$M_2$ is normal, so that~$B$ fits into exactly one extension of the type
\[M \hookrightarrow B \twoheadrightarrow Q\]
with $\dm M = 8$ and $\dm Q = 4$. In this case, $M = M_ 2 \cong K^{D_8}$,
while \mbox{$Q \cong K[\Z_2 \times \Z_2]$}. This extension is abelian, but neither central nor cocentral.

In Paragraph~\ref{normal4}, we determine the four-dimensional Hopf subalgebras of~$B$. There are seven, but only one of them is normal. Consequently, the Hopf algebra~$B$ fits into exactly one extension of the type
\[P \hookrightarrow B \twoheadrightarrow F \]
with $\dm P = 4$ and $\dm F = 8$. In this case, $P \cong K[\Z_2\times \Z_2]$ and $F \cong K[D_8]$. This extension is abelian, but neither central nor cocentral.

Besides the Yetter-Drinfel'd Hopf algebra~$A$, the authors considered in~\cite{KaSo2}, Section~3 a second eight-dimensional Yetter-Drinfel'd Hopf algebra over \mbox{$K[\Z_2 \times \Z_2]$}, denoted here by~$A'$. These two algebras are certainly not isomorphic, because~$A$ is commutative, while~$A'$ is not. However, the corresponding Radford biproducts~$B$ and~$B'$ are isomorphic, as we show in Section~\ref{Sec:BiprodSec}, so that it is not necessary to carry out the same analysis for~$B'$.

Let us now state the conventions that are used throughout this article. Our base field~$K$ will be algebraically closed of characteristic zero. The multiplicative group of~$K$ will be denoted by~$K^\times = K \setminus \{0\}$.

All vector spaces will be defined over~$K$. The dual space of a vector space~$V$ will be denoted by~$V^* = \Hom(V,K)$, and the dual map of a linear map~$f$, also called its transpose, will be denoted by~$f^*$.

The group algebra of a group~$G$, which we will also call its group ring, will be denoted by~$K[G]$. The corresponding dual space~$K[G]^*$, the dual group ring, will be denoted by~$K^G$. The character group of a group~$G$ will be denoted by \mbox{$\hat{G} = \Hom(G,K^\times)$}. Its elements will be interchangeably called multiplicative characters, one-dimensional characters, or one-dimensional representations, a terminology that we will use not only for groups, but also for algebras. Clearly, the multiplicative characters of the group algebra correspond bijectively to the multiplicative characters of the group via restriction.

All algebras are assumed to have a unit element, and algebra homomorphisms are assumed to preserve these unit elements. The center of an algebra~$A$ will be denoted by~$Z(A)$, which should not be confused with a certain quotient Hopf algebra~$Z$ already mentioned above. The set of group-like elements in a coalgebra~$A$ will be denoted by~${\bf G}(A)$. The subalgebra generated by elements~$a_1,\ldots,a_n$ of an algebra~$A$ will be denoted by~$K \langle a_1,\ldots,a_n \rangle$, whereas the subgroup of a group~$G$ generated by elements~$g_1,\ldots,g_n$  will be denoted by~$\langle g_1,\ldots,g_n \rangle$. The augmentation ideal of a Hopf algebra~$H$ will be denoted by~$H^+ \deq \ker(\eh)$. The image of an element~$a$ in a quotient space will be denoted by~$\bar{a}$. Finally, the symbol~$\subset$ denotes non-strict inclusion.

As already discussed, the article is divided into sections, which are divided further into relatively small paragraphs. A reference to Proposition~2.2 refers to the unique proposition in Paragraph~2.2, and definitions, theorems, lemmas, and corollaries are referenced in the same way.

The material discussed here and in our previous article~\cite{KaSo2} was presented at the AMS Spring Southeastern Sectional Meeting in Auburn in March~2019, in two consecutive talks. It was also presented at the International Workshop on Hopf Algebras and Tensor Categories in Nanjing in September~2019. The travel of the second author to these conferences as well as his work on this article were supported by NSERC grant RGPIN-2017-06543. The work of the first author on this article was supported by a Faculty Summer Research Grant from the College of Science and Health at DePaul University.

\section{Preliminaries} \label{Sec:Prelim}
\subsection[Yetter-Drinfel'd category]{} \label{YetDrinf}
In this section, we consider a finite-dimensional Hopf algebra~$H$ with coproduct~$\dhh$, counit~$\eh$, and antipode~$\sh$. Recall that a Yetter-Drinfel'd module over~$H$ is a vector space~$V$ that is both a left $H$-module and a left $H$-comodule in such a way that the compatibility condition
\begin{equation*}
h_{\1} v^{\1} \ot h_{\2}.v^{\2} =
\left(h_{\1}.v \right)^{\1} h_{\2} \ot \left( h_{\1}.v \right)^{\2}
\end{equation*}
is satisfied for all $h\in H$ and~$v\in V$, where we have used the notation
\begin{equation*}
\delta(v) = v^{\1} \ot v^{\2} \in H \ot V
\end{equation*}
for the coaction~$\delta$. The category $_{H}^{H}YD$ of all Yetter-Drinfel'd modules, whose morphisms are the $H$-linear and colinear maps, is a braided monoidal category with the braiding
$\sigma_{V,W} \colon V \ot W \rightarrow W \ot V$ given by
\begin{equation*}
\sigma_{V,W} \left( v \ot w \right) = v^{\1}.w \ot v^{\2}
\end{equation*}
We collect a number of facts about Yetter-Drinfel'd modules that will be needed in the sequel:
\begin{parlist}
\item
The dual space~$V^*$ of a finite-dimensional Yetter-Drinfel'd module~$V$ can be considered as a Yetter-Drinfel'd module over the dual Hopf algebra~$H^*$ by using the module structure
\begin{equation}
\gamma.\varphi = \left( \gamma \ot \varphi \right) \circ \delta  \label{dualaction}
\end{equation}
for $\gamma \in H^*$ and $\varphi \in V^*$, and the comodule structure 
$\delta(\varphi) = \varphi^{\1} \ot \varphi^{\2} \in H^* \ot V^*$ with the defining property 
\begin{equation}
\varphi^{\1}(h) \; \varphi^{\2}(v) = \varphi(h.v) \label{dualcoaction}
\end{equation}
where $\varphi \in V^*$, $h \in H$, and~$v \in V$ (cf.~\cite{SoYp}, Proposition~1.2 and Lemma~1.3). With respect to these structures, the transpose of an $H$-linear and $H$-colinear map is $H^*$-linear and $H^*$-colinear.

\item
If~$A$ and~$A'$ are algebras in the Yetter-Drinfel'd category, so that their multiplication maps~$\mu_{A}$ and~$\mu_{A'}$ as well as their unit maps~$\eta_{A}$ and $\eta_{A'}$ are morphisms in the category, then their tensor product~$A \ot A'$ is also an algebra in this category when endowed with the multiplication map
\[\mu_{A \ot A'} = \left( \mu_{A} \ot \mu_{A'} \right) \circ \left(\id_{A} \ot \, \sigma_{A',A} \ot \id_{A'} \right)\]
and the unit map $\eta_{A \ot A'} = \eta_{A} \ot \eta_{A'}$.

\item
A Hopf algebra in the Yetter-Drinfel'd category~$_{H}^{H}YD$, or a Yetter-Drinfel'd Hopf algebra for short, is an algebra~$A$ that is simultaneously a coalgebra, both in this category, so that comultiplication and counit are multiplicative. Here, the algebra structure on~$A \ot A$ is the one introduced above, so that the multiplicativity of the coproduct means that
\[\Delta_{A} \circ \mu_{A} =
\left( \mu_{A} \ot \mu_{A}\right) \circ \left(\id_{A} \ot \, \sigma_{A,A} \ot \id_{A} \right) \circ
\left( \Delta _{A} \ot \Delta _{A} \right)\]
or, if we use the notation $\Delta _{A}\left( a \right) = a_{\1} \ot a_{\2}$ for the coproduct,
\[\Delta _{A}\left( a a' \right) =
a_{\1} \left( a_{\2}{}^{\1}.a'_{\1}\right) \ot a_{\2}{}^{\2} a'_{\2}\]
We also require that~$A$ has an antipode~$S_{A}$, i.e., a convolution inverse of the identity. As shown in~\cite{KaSo2}, Lemma~1.2, such an antipode is automatically a morphism in the category.
\end{parlist}

\subsection[Biproducts]{}  \label{Biprod}
Yetter-Drinfel'd Hopf algebras give rise to ordinary Hopf algebras via the Radford biproduct construction, which was originally defined in \cite{RadfProj}, Theorem~1:
\begin{defn}
Let $A$ be a Yetter-Drinfel'd Hopf algebra over~$H$. The Radford biproduct
$B \deq A\star H$ is the tensor product of~$A$ and~$H$ as a vector space, endowed with the smash product multiplication
\begin{eqnarray} \label{Multip}
\left( a \star h\right) \left( a' \star  h' \right) \deq
a \left( h_{\1}.a'\right) \star h_{\2} h'
\end{eqnarray}
with corresponding unit element~$1_{B} \deq 1_{A} \star 1_{H}$, cosmash product comultiplication
\begin{eqnarray} \label{Comult}
\db \left( a \star h \right) \deq
\left( a_{\1} \star a_{\2}{}^\1 h_{\1} \right) \ot \left(a_{\2}{}^\2 \star h_{\2} \right)
\end{eqnarray}
with corresponding counit
$\eb (a \star h) \deq \ea (a) \eh(h)$, and antipode
\begin{eqnarray} \label{Antip}
\sbb \left( a \star h \right) \deq
\left( \A \star \sh(a^{\1} h)
\right) \left( \sa (a^{\2}) \star \HH \right)
\end{eqnarray}
With respect to these structures, it is an ordinary Hopf algebra. The notation~\mbox{$a \star h$} instead of~$a \ot h$ is used to signify that the algebra and coalgebra structures are not the canonical structures on the tensor product.
\end{defn}

Two results about Radford biproducts will be needed in the sequel:
\begin{parlist}
\item
If $A$ is a finite-dimensional Yetter-Drinfel'd Hopf algebra over~$H$, then~$A^*$ is a Yetter-Drinfel'd Hopf algebra over~$H^*$ with respect to the module and comodule structures described in Paragraph~\ref{YetDrinf} and the standard dualizations of the algebra and coalgebra structures (cf.~\cite{SoYp}, loc.~cit.), so that we can form the biproduct $A^* \star H^*$. It is isomorphic to~$B^*$ via the canonical isomorphism of the underlying vector spaces.

\item
The biproduct~$B$ is semisimple if~$A$ and~$H$ are semisimple. Similarly, it is cosemisimple if~$A$ and~$H$ are cosemisimple (cf.~\cite{RadfProj}, Proposition~3 and Proposition~4; \cite{SoRib}, Paragraph~2.14 and the references cited there).
\end{parlist}

Further details on Yetter-Drinfel'd Hopf algebras and biproducts can be found in~\cite{M}, \S~10.6, \cite{RadfProj}, and~\cite{SoYp}.

\subsection[The universal property of a crossed product]{} \label{UnivPropCross}
A generalization of the smash products just encountered in Paragraph~\ref{Biprod} are crossed products (cf.~\cite{M}, Definition~7.1.1): Here we have, for a Hopf algebra~$H$ and an algebra~$A$, not only a weak action, or measuring, $\rightharpoonup \colon H \ot A \to A$, but also a cocycle~\mbox{$\tau \colon H \times H \to A$}, and the multiplication on the underlying vector space~$A \ot H$, which is then commonly denoted by $A \#_{\tau } H$, is given by the formula
\[(a \# h)(a' \# h') \deq a (h_\1 \rightharpoonup a') \tau(h_\2, h'_\1) \# h_\3 h'_\2\]
Such a crossed product has the following universal property:
\begin{prop}
Suppose that~$B$ is an algebra and that~$f_A \colon A \ra B$ as well as \mbox{$f_H \colon H \ra B$} are $K$-linear maps. Then the following assertions are equivalent:
\begin{enumerate}
\item
There exists an algebra homomorphism $f \colon A \#_{\tau } H \ra B$ such that
$$f(a \# \HH) = f_A(a) \qquad \text{and} \qquad f(\A \# h) = f_H(h)$$

\item
$f_A$ is an algebra homomorphism, $f_H(\HH) = \B$, and we have
$$f_A(h_\1 \rightharpoonup a) f_H(h_\2) = f_H(h) f_A(a)$$
as well as
$$f_A(\tau(h_\1, h'_\1)) f_H(h_\2 h'_\2) = f_H(h) f_H(h')$$
\end{enumerate}
\end{prop}
\begin{pf}
To show that the first assertion implies the second, we define the inclusion maps $i_A \colon A\ra A \#_{\tau } H$ and $i_H \colon H \ra  A \#_{\tau } H$ 
via~$i_A(a) \deq a \# \HH$ and \mbox{$i_H(h) \deq \A \# h$}, and also define $f_A \deq f \circ i_A$ and $f_H \deq f \circ i_H$. Then we have $f_H(\HH) = f(\A \# \HH) = \B$ and
\begin{align*}
f_H(h) f_A(a)  &= f(\A \# h)f(a \# \HH) =
f\left( (\A \# h)(a \# \HH)\right)\\
& = f\left( (h_\1 \rightharpoonup a) \# h_\2 \right) = f_A(h_\1 \rightharpoonup a) f_H(h_\2)
\end{align*}
as well as
\begin{align*}
f_H(h) f_H(h')  &=
f(\A \# h) f(\A \# h')=f\left( (\A \# h) (\A \# h')\right)  \\
&=
f \left( \tau(h_\1, h'_\1) \# h_\2 h'_\2\right) = f_A(\tau(h_\1, h'_\1)) f_H(h_\2 h'_\2)
\end{align*}

To show that the second assertion implies the first, we define
\[f(a \# h) \deq f_A(a) f_H(h)\]
Then~$f$ preserves the unit element, and we have
\begin{align*}
f((a \# h)(a' \# h')) &= f(a (h_\1 \rightharpoonup a') \tau(h_\2, h'_\1) \# h_\3 h'_\2) \\
&= f_A(a) f_A(h_\1 \rightharpoonup a') f_A(\tau(h_\2, h'_\1)) f_H(h_\3 h'_\2) \\
&= f_A(a) f_A(h_\1 \rightharpoonup a') f_H(h_\2) f_H(h') \\
&= f_A(a) f_H(h) f_A(a') f_H(h') = f(a \# h) f(a' \# h')
\end{align*}
as asserted.
\qed
\end{pf}

There are two observations about this universal property that are worth recording. First, if the equation
\[f_A(h_\1 \rightharpoonup a) f_H(h_\2) = f_H(h) f_A(a)\]
that appears in the universal property is satisfied for all~$h \in H$, but only for a
specific~$a \in A$, and in the same way for another specific element~$a' \in A$, then it holds for their product~$aa'$, as we have by the measuring property that
\begin{align*}
f_A(&h_\1 \rightharpoonup (aa')) f_H(h_\2) =
f_A((h_\1 \rightharpoonup a) (h_\2 \rightharpoonup a')) f_H(h_\3) \\
&= f_A(h_\1 \rightharpoonup a) f_A(h_\2 \rightharpoonup a') f_H(h_\3) \\
&= f_A(h_\1 \rightharpoonup a) f_H(h_\2) f_A(a') = f_H(h) f_A(a) f_A(a') = f_H(h) f_A(aa')
\end{align*}
Arguing inductively, this shows that it is sufficient to verify this equation for all~$a$ in a generating set of~$A$, as long as it holds for all~$h \in H$ for each of these generators.

For the second observation, note that~$i_H$ is convolution-invertible (cf.~\cite{M}, Proposition~7.2.7). If $\bar{i}_H$ denotes its convolution-inverse, then
$\bar{f}_H \deq f \circ \bar{i}_H$ is the convolution-inverse of
$f_H = f \circ i_H$. Thus  the two main conditions of the second assertion can also be stated in the form
$$f_A(h \rightharpoonup a) = f_H(h_\1) f_A(a) \bar{f}_H(h_\2)$$
and
$$f_A(\tau(h, h')) = f_H(h_\1) f_H(h'_\1) \bar{f}_H(h_\2 h'_\2)$$

In the case that the cocycle is trivial in the sense that $\tau(h, h') = \eh(h) \eh(h') \A$, the crossed product reduces to the smash product. Then the twisted module condition (cf.~\cite{M}, Lemma~7.1.2) reduces to the usual module condition, so that the weak action is an ordinary action, which we will denote by a dot again. The preceding proposition now yields the following universal property of the smash product:
\begin{cor}
Suppose that~$B$ is an algebra and that~$f_A \colon A \ra B$ and \mbox{$f_H \colon H \ra B$} are algebra homomorphisms such that
$$f_A(h_\1.a) f_H(h_\2) = f_H(h) f_A(a)$$
Then the map $f \colon A \# H \ra B$ defined via
$$f(a \# h) \deq f_A(a) f_H(h)$$
is an algebra homomorphism.
\end{cor}

For the smash product, there is a second point to be made, which is similar, but not identical, to the one made above in the case of a general crossed product: If the hypothesis in the preceding corollary holds for all~$a \in A$ and two elements~$h, h' \in H$, it holds for their product, because the weak action is then an ordinary action and we therefore have
\begin{align*}
f_A(h_\1 h'_\1.a) &f_H(h_\2 h'_\2) = f_A(h_\1 h'_\1.a) f_H(h_\2) f_H(h'_\2) \\
&= f_H(h) f_A(h'_\1.a) f_H(h'_\2) \\
&= f_H(h) f_H(h') f_A(a) = f_H(h h') f_A(a)
\end{align*}
Arguing inductively, this shows that it is sufficient to verify this hypothesis for all~$h$ in a generating set of~$H$, as long as it holds for all~$a \in A$ for each of these generators.

\subsection[Extensions]{} \label{Extens}
As we will explain in Paragraph~\ref{Bicross} below, crossed products appear naturally when studying Hopf algebra extensions, which are defined as follows (cf.~\cite{MaQu}, Definition~5.6):
\begin{defn}
Consider a sequence
\begin{equation*}
N \overset{\iota_N}{\hookrightarrow} B \overset{\pi_F}{\twoheadrightarrow} F
\end{equation*}
of finite-dimensional Hopf algebras and Hopf algebra maps~$\iota_N$ and~$\pi_F$, where~$\iota_N$ is injective and~$\pi_F$ is surjective. This sequence is called an extension of~$F$ by~$N$ if
\[\iota_N(N) =
\left\{ b\in B \mid b_{\1} \ot \pi_F(b_{\2}) = b \ot 1_{F}\right\} \]
In this situation, $\iota_N(N)$ is a normal Hopf subalgebra of~$B$ (cf.~\cite{AndNot}, Section~1.2).
\end{defn}

Hopf algebra extensions can have additional properties:
\begin{parlist1}
\item
The extension is called abelian if~$N$ is commutative and~$F$ is
cocommutative. Since~$K$ is algebraically closed of characteristic zero, it follows from~\cite{M}, Theorem~2.3.1 that~$N \cong K^\Gamma$ and $F \cong K[L]$ for two finite groups~$\Gamma$ and~$L$.

\item
The extension is called central if $\iota_N(N) \subset Z(B)$, the center of~$B$. This happens in particular if~$N$ has dimension~$2$, since~$N$ is then necessarily semisimple and is spanned by its unit element and its integral, both of which are central in view of~\cite{Ma0}, Lemma~2.16.

\item
The extension is called cocentral if $\pi_F^*(F^*) \subset Z(B^*)$, the center of~$B^*$.

\item
We say that two extensions $N \hookrightarrow B \twoheadrightarrow F$ and $N \hookrightarrow B' \twoheadrightarrow F$ are equivalent if there is a Hopf algebra isomorphism $f \colon B\rightarrow B'$ that induces the identity maps on~$N$ and~$F$ (cf.~\cite{Ma4}, Definition~1.4).
\end{parlist1}

Further details on Hopf algebra extensions can be found in~\cite{AD}, \cite{AndNot}, \cite{HofstDiss}, \cite{Hofst}, \cite{MaQu}, and~\cite{Ma4}.

\subsection[Bicrossed products]{} \label{Bicross}
The connection between crossed products and Hopf algebra extensions arises via bicrossed products, which are defined as follows:
\begin{defn}
Let $N$ and $F$ be Hopf algebras equipped with a weak action, or measuring,
$\rightharpoonup \colon F\ot N\rightarrow N$, a cocycle $\tau \colon F \times F \rightarrow N$, and dually a weak coaction $\rho \colon F \rightarrow F \ot N$ and a dual cocycle
$\kappa \colon F \rightarrow N \ot N$. The bicrossed product $N \#_{\tau}^{\kappa} F$ is a Hopf algebra with~$N \ot F$ as underlying vector space, crossed product multiplication
\[(x \# t) (y \# s) \deq x (t_{\1} \rightharpoonup y) \tau(t_{\2}, s_{\1}) \# t_{\3} s_{\2} \]
as in Paragraph~\ref{UnivPropCross}, and dually cocrossed product comultiplication
\[\Delta (x\# t) \deq
x_{\1} \kappa_1(t_{\1}) \# t_{\2}{}^{\1}
\ot x_{\2} \kappa_2(t_{\1}) t_{\2}{}^{\2} \# t_{\3}\]
where we have, in analogy to Paragraph~\ref{YetDrinf}, used the notation
\[\rho(t) = t^{\1} \ot t^{\2} \in F \ot N \]
even though this is a right, not a left, weak coaction. For the dual cocycle, we have used the notation
\[\kappa(t) = \kappa_{1}(t) \ot \kappa_{2}(t) \in N \ot N\]
\end{defn}

The necessary compatibility conditions on the four structure elements $\rightharpoonup$, $\tau$, $\rho$, and~$\kappa$ can be found in~\cite{AD}, Theorem~2.20; \cite{AndNot}, Section~3.1; or~\cite{HofstDiss}, Kapitel~5.

Every bicrossed product yields the extension
\begin{equation*}
N \overset{\iota_N}{\hookrightarrow} N \#_{\tau}^{\kappa} F
\overset{\pi_F}{\twoheadrightarrow} F
\end{equation*}
where $\iota_N(x) \deq x \# 1_F$ and~$\pi_F(x \# t) \deq \varepsilon_N(x) t$. In the finite-dimensional case considered here, in fact every extension is equivalent to such a bicrossed product extension (cf.~\cite{AD}, Theorem~3.2.14; \cite{AndNot}, Proposition~3.1.12 and Theorem~3.1.17; \cite{MaQu}, Page~130).

It is not difficult to see from the formula for the crossed product multiplication given above that, in a central bicrossed product extension, the weak action of~$F$ on~$N$ must be trivial. Dually, for a cocentral extension $B = N \#_{\tau}^{\kappa} F$,
the weak coaction~$\rho$ is trivial, because in~$B^*$ the weak action of~$N^*$ on~$F^*$ is trivial.

\section{The biproduct in the first case} \label{Sec:BiprodFirst}
\subsection[The Yetter-Drinfel'd Hopf algebra in the first case]{} \label{YetDrinfFirst}
If we use the notation~$\Z_2 = \{0,1\}$ for the group with two elements, the group
$G \deq \Z_2 \times \Z_2$ contains the four elements
\[g_1 \deq (0,0) \qquad g_2 \deq (1,0) \qquad g_3 \deq (0,1) \qquad  g_4 \deq (1,1)\]
with $g_1 = 1_G$ as its unit element when we write~$G$ multiplicatively, as we will in the sequel.
In \cite{KaSo2}, Section~2, the authors defined a commutative Yetter-Drinfel'd Hopf algebra~$A$ over the group algebra~$H \deq K[G]$ as the algebra generated by two commuting variables~$x$ and~$y$ subject to the defining relations
\[x^4 = 1 \qquad \qquad y^2 = \frac{1}{2}(1 + \zeta x + x^2 - \zeta x^3)\]
where~$\zeta$ is a not necessarily primitive fourth root of unity. With the help of a primitive fourth root of unity~$\iota$, it is possible to introduce the elements $\omega_1 \deq \ua$,
\[\omega_2 \deq \frac{\K}{2}(1 + \iota \zeta^2) x + \frac{\K}{2}(1 - \iota \zeta^2) x^3 \qquad
\omega_3 \deq \frac{\K}{2}(1 - \iota \zeta^2) x + \frac{\K}{2}(1 + \iota \zeta^2) x^3\]
$\omega_4 \deq x^2$, and
\[\eta_1 \deq y \qquad \eta_2 \deq x^3 y \qquad \eta_3 \deq x^2 y \qquad \eta_4 \deq x y\]
According to \cite{KaSo2}, Proposition~2.4, these eight elements form a basis of~$A$, and the coalgebra structure is determined by the fact that they are group-like. As also shown there, the set $\{\omega_1 , \omega_2, \omega_3, \omega_4\}$ is closed under multiplication and in fact forms a group isomorphic
to~$\Z_2 \times \Z_2$. The products of all basis elements are computed there; here we only record that $\omega_4 \eta_1 = \eta_3$ and $\omega_4 \eta_2 = \eta_4$. In~\cite{KaSo2}, Paragraph 2.2, the 
$H$-action on~$A$ is introduced on the generators as
\[g_2.x = x^3 \qquad \qquad g_2.y = x^3 y \qquad \qquad
g_3.x = x \qquad \qquad g_3.y = x^2 y\]
which according to~\cite{KaSo2}, Lemma~2.4 means for the basis elements that
\[g_3.\omega_i = \omega_i \qquad \qquad  g_2.\omega_4 = \omega_4 \qquad \qquad g_2.\omega_2 = \omega_3 \qquad \qquad  g_i.\eta_{g_j} = \eta_{g_i g_j}\]
where we have used the notation~$\eta_{g_i} \deq \eta_i$. In particular, we have
$g_3.\eta_i = \omega_4 \eta_i$ for all $i=1, 2, 3, 4$. The coaction of~$H$ on~$A$ is derived from the action by defining
\begin{equation}
\delta(a) \deq \frac{1}{4} \sum_{g,g' \in G} \theta(g,g') \, g \ot g'.a \label{coaction}
\end{equation}
where~$\theta$ is the nondegenerate symmetric bicharacter described in~\cite{KaSo2}, Paragraph~2.2 and determined by the condition
\[\begin{pmatrix}
\theta(g_2,g_2) & \theta(g_2,g_3) \\
\theta(g_3,g_2) & \theta(g_3,g_3)
\end{pmatrix} =
\begin{pmatrix}
\zeta^2 & -1 \\
-1 & 1
\end{pmatrix} \]
on its fundamental matrix.

We will need the following two explicit forms, the first of which is recorded in~\cite{KaSo2}, Paragraph~4.3:
\begin{align*}
\delta(\omega_2) &= \frac{1}{2} (g_1 + g_3) \ot \omega_2 
+ \frac{1}{2} (g_1 - g_3) \ot \omega_3\\
\delta(\eta_1) &= \frac{1}{4}(g_1+g_2+g_3+g_4) \ot \eta_1 
+ \frac{1}{4}(g_1+\zeta^2 g_2-g_3-\zeta^2 g_4) \ot \eta_2\\
&\quad + \frac{1}{4}(g_1-g_2+g_3-g_4) \ot \eta_3 + \frac{1}{4}(g_1-\zeta^2 g_2-g_3+\zeta^2 g_4)\ot  \eta_4
\end{align*}

Because~$\eta_1 = y$, the last formula also yields the value of the coaction on one of the generators. According to~\cite{KaSo2}, Proposition~2.2, the value on the other generator is
\begin{align*}
&\delta(x) = \frac{\K}{2} (g_1 + g_3) \ot x + \frac{\K}{2} (g_1 - g_3) \ot x^3 \\
&\delta(x^3) = \frac{\K}{2} (g_1 - g_3) \ot x + \frac{\K}{2} (g_1 + g_3) \ot x^3
\end{align*}
where it is also stated that~$x^2 = \omega_4$ is coinvariant.

It is possible to expand the powers of~$x$ explicitly in terms of the group-like elements: If we multiply the definition of~$\omega_2$ by $1 + \iota \zeta^2$ and the definition of~$\omega_3$ by $1 - \iota \zeta^2$, we obtain
\[(1 + \iota \zeta^2) \omega_2 = x^3 + \iota \zeta^2 x \qquad \qquad
(1 - \iota \zeta^2) \omega_3 = x^3 - \iota \zeta^2 x\]
Adding these two formulas, we get
\begin{align*}
x^3 = \frac{1}{2} (1 + \iota \zeta^2) \omega_2 + \frac{1}{2} (1 - \iota \zeta^2) \omega_3
\end{align*}
Subtracting these formulas, we obtain
$2 \iota \zeta^2 x = (1 + \iota \zeta^2) \omega_2 - (1 - \iota \zeta^2) \omega_3$ or
\begin{align*}
x = \frac{1}{2} (1 - \iota \zeta^2) \omega_2 + \frac{1}{2} (1 + \iota \zeta^2) \omega_3 
\end{align*}
It follows from these considerations that the space
\[C \deq \Span(1, x, x^2 ,x^3) = \Span(\omega_1, \omega_2, \omega_3, \omega_4)\]
is a Yetter-Drinfel'd Hopf subalgebra of~$A$. As shown in~\cite{KaSo2}, Paragraph~4.3, it is the so-called core of~$\eta_1$; it is trivial, but not completely trivial in the sense explained there.

\subsection[The algebra structure of the biproduct]{} \label{AlgBiprod}
\enlargethispage{8pt}
For this specific Yetter-Drinfel'd Hopf algebra~$A$ over the group algebra~$H$, we now form the biproduct~$B \deq A \star H$. It is generated by the four elements
\[u \deq x \star \HH \qquad \quad v \deq y \star \HH \qquad \quad
r \deq \A \star g_2 \qquad \quad s \deq \A \star g_3\]
These generators can be used to give a presentation of~$B$:
\begin{prop}
The four generators satisfy the relations
\begin{proplist}
\item
$\displaystyle
u^4 = 1, \qquad uv = vu, \qquad v^2 = \frac{1}{2}(1 + \zeta u + u^2 - \zeta u^3)$

\item
$\displaystyle
r^2 = 1, \qquad rs = sr, \qquad s^2 = 1$

\item
$\displaystyle
r u = u^3 r, \qquad r v = u^3 v r, \qquad s u = u s, \qquad s v = u^2 v s$
\end{proplist}
These relations are defining.
\end{prop}
\begin{pf}
The first set of relations is a direct consequence of the defining relations of~$A$, and the second set is a consequence of the relations inside the group~$G$. The third set follows from the definition of the action of~$G$ on~$A$ and the multiplication in a smash product from Equation~(\ref{Multip}): We have
\[ru = (\A \star g_2) (x \star \HH) = (g_2.x) \star g_2 = x^3 \star g_2 = u^3 r \]
and similarly $r v = (g_2.y) \star g_2 = (x^3 y) \star g_2$ as well as
$s u = (g_3.x) \star g_3 = u s$ and $s v = (g_3.y) \star g_3 = (x^2 y) \star g_3$.
To see that the relations are defining, we consider the abstract algebra given by this presentation. The computations just carried out yield a surjective algebra homomorphism from this abstract algebra to~$B$, which has dimension~$32$. On the other hand, because every pair of the abstract generators satisfies a commutation relation, the abstract algebra is spanned by the 32~elements~$u^i v^j r^k s^l$, where the exponents take the values $i=0,1,2,3$ and~$j,k,l=0,1$. Therefore the abstract algebra has at most dimension~$32$. As it maps surjectively onto a space of dimension~$32$, it must have exactly this dimension, and the map must be an isomorphism.
\qed
\end{pf}

The basis elements of~$A$ can be used in a similar way to define
\[c_i \deq \omega_i \star \HH \qquad  \qquad  d_i \deq \eta _i \star \HH \qquad \qquad
h_i \deq 1_A \star g_i\]
for $i=1,2,3,4$. Then the elements $c_i h_j$ together with the elements $d_i h_j$,
for~$i, j = 1,2,3, 4$, form a basis of~$B$. Note that $c_1 = h_1 = \B$. Except for~$u$, the above generators are among these elements, as~$v=d_1$, $r=h_2$, and~$s=h_3$. Because the canonical map from~$A$ to~$B$ is an algebra homomorphism, the relations that we have derived in Paragraph~\ref{YetDrinfFirst} yield corresponding relations in~$B$. In particular, the elements~$c_i$ commute with elements~$d_j$, and we have
\[c_3 = c_4 c_2 \qquad \qquad d_3 = c_4 d_1 \qquad \qquad d_4 = c_4 d_2  \]
As in the proof of the preceding proposition, Equation~(\ref{Multip}), the formula for the multiplication in a smash product, yields the commutation relations
\[h_3 c_i = c_i h_3 \qquad \quad  h_j c_4  = c_4 h_j \qquad \quad
h_2 c_2 = c_3 h_2 \qquad \quad  h_id_{g_j} = d_{g_i g_j}h_i\]
where we have used the analogous notation~$d_{g_i} \deq d_i$. We also note that
\begin{align} \label{c2h2}
(c_2 h_2)^2 = c_2 h_2 c_2 h_2 = c_2 c_3 h_2 h_2 = c_2 c_3 = c_4
\end{align}
so that $(c_2 h_2)^3 = c_3 h_2$ and $(c_2 h_2)^4 = 1$.

\subsection[The coalgebra structure of the biproduct]{} \label{CoalgBiprod}
The form of the coproduct can be derived from Equation~(\ref{Comult}). According to~\cite{RadfProj}, Fact~2.11, the group-like elements of~$B$ are of the form $a \star g_j$, where $a \in A$ is group-like and coinvariant. This implies that we have $\Delta_B(c_4) = c_4 \ot c_4$ and $\Delta _B(h_j) = h_j \ot h_j$. The formula for the coaction of~$\omega_2$ yields 
\begin{align*}
\db(c_2) &= \frac{1}{2} c_2(h_1 + h_3) \ot c_2 + \frac{1}{2} c_2(h_1 - h_3) \ot c_3 \\
&= \frac{1}{2}\left( (h_1 + h_3) \ot c_1 + (h_1 - h_3) \ot c_4 \right)\left( c_2 \ot c_2\right)
\end{align*}
while the formula for the coaction of~$\eta_1$ yields
\begin{align*}
\db(d_1) &= \frac{1}{4}d_1(h_1 + h_2 + h_3 + h_4) \ot d_1
+ \frac{1}{4} d_1 (h_1 + \zeta^2 h_2 - h_3 - \zeta^2 h_4) \ot d_2 \\
&\quad + \frac{1}{4} d_1 (h_1 - h_2 + h_3 - h_4) \ot d_3
+ \frac{1}{4} d_1 (h_1 - \zeta^2 h_2 - h_3 + \zeta^2 h_4) \ot d_4 \\
&= \frac{1}{4} \left( d_1 \ot d_1\right) \left((h_1+h_3) \ot (c_1 + c_4)
+ (h_2+h_4) \ot  (c_1 - c_4) \right)  \\
&\quad +\frac{1}{4} \left( d_1 \ot d_2 \right)
\left((h_1 - h_3) \ot (c_1 + c_4) + \zeta^2(h_2 - h_4) \ot (c_1 - c_4)  \right)
\end{align*}
It should be noted that the coproduct of~$B$ is completely determined by the formulas above, since they also allow to compute the coproduct of~$c_3 = c_4 c_2$ and~$d_j = h_j d_1 h_j$ by multiplicativity. The counit is given on these elements by~$\eb(c_i) = \eb(d_i) = \eb(h_i) = 1$ for all $i=1,2,3,4$.

An alternative way of determining the coproduct is to record its values on the generators:
\begin{prop}
We have $\db(r) = r \ot r$ and $\db(s) = s \ot s$. Furthermore, we have
\begin{align*}
\db(u) &= \frac{\K}{2} (u \ot u + u \ot u^3 +  u^3 s \ot u - u^3 s \ot u^3)
\end{align*}
and
\begin{align*}
\db(v) &= \frac{1}{4} v (1 + r + s + rs) \ot v
+ \frac{1}{4} v (1 - \zeta^2 r - s + \zeta^2 rs) \ot u v \\
&\quad + \frac{1}{4} v (1 - r + s - rs) \ot u^2 v
+ \frac{1}{4} v (1 + \zeta^2 r - s - \zeta^2 rs) \ot u^3 v
\end{align*}
The counit is given on generators by~$\eb(u) = \eb(v) = \eb(r) = \eb(s) = 1$.
\end{prop}
\begin{pf}
We have already said that $r=h_2$ and $s=h_3$ are group-like. Since $v=d_1$, the formula for~$\db(v)$ is just a restatement of the corresponding one for~$\db(d_1)$. According to~\cite{KaSo2}, Proposition~2.3, we have
\[\da(x) =  \frac{\K}{2} (x \ot x + x \ot x^3 + x^3 \ot x - x^3 \ot x^3)\]
With the help of the formulas for~$\delta(x)$ and~$\delta(x^3)$ recalled in Paragraph~\ref{YetDrinfFirst}, Equation~(\ref{Comult}) then yields
\begin{align*}
\db(u) = {} &
\frac{\K}{4} (x \star (g_1 + g_3)) \ot (x \star \HH)
+ \frac{\K}{4} (x \star (g_1 - g_3)) \ot (x^3 \star \HH)\\
&+ \frac{\K}{4} (x \star (g_1 - g_3)) \ot (x \star \HH)
+ \frac{\K}{4} (x \star (g_1 + g_3)) \ot (x^3 \star \HH)\\
&+ \frac{\K}{4} (x^3 \star (g_1 + g_3)) \ot (x \star \HH)
+ \frac{\K}{4} (x^3 \star (g_1 - g_3)) \ot (x^3 \star \HH)\\
&- \frac{\K}{4} (x^3 \star (g_1 - g_3)) \ot (x \star \HH)
- \frac{\K}{4} (x^3 \star (g_1 + g_3)) \ot (x^3 \star \HH)
\end{align*}
By adding the first and the second as well as the third and the fourth line, this expression becomes
\begin{align*}
\db(u) = {} &
\frac{\K}{2} (x \star g_1) \ot (x \star \HH)
+ \frac{\K}{2} (x \star g_1) \ot (x^3 \star \HH) \\
&+ \frac{\K}{2} (x^3 \star g_3) \ot (x \star \HH)
- \frac{\K}{2} (x^3 \star g_3) \ot (x^3 \star \HH)
\end{align*}
Using that~$g_1 = \HH$, we get the formula for~$\db(u)$. The form of the counit follows directly from the defining equations~$\ea(x) = \ea(y) = 1$ in~\cite{KaSo2}, Paragraph~2.3.
\qed
\end{pf}

\subsection[The antipode of the biproduct]{} \label{AntBiprod}
To complete the description of~$B$, we compute the values of the antipode on the generators:
\begin{prop}
We have $\sbb(u) = \frac{\K}{2} (u + s u + u^3 - s u^3)$ and
\begin{align*}
\sbb \left( v \right) = {} &
\frac{1}{4} \left(1 + u + u^2 + u^3 \right) v
+ \frac{1}{4} \left(\zeta 1 + \zeta^3 u - \zeta u^2 - \zeta^3 u^3 \right) r v \\
&+ \frac{1}{4} \left(1 - u + u^2 - u^3 \right) s v
+ \frac{1}{4}  \left(- \zeta 1 + \zeta^3 u
+ \zeta u^2 - \zeta^3 u^3 \right) r s v
\end{align*}
as well as $\sbb(r) = r$ and~$\sbb(s) = s$.
\end{prop}
\begin{pf}
The forms of~$\sbb(r)$, $\sbb(s)$, and~$\sbb(u)$ follow from Equation~(\ref{Antip}), where the third also requires the equation for~$\delta(x)$ given in Paragraph~\ref{YetDrinfFirst} as well as the facts that~$\sa(x) = x$ and $\sa(x^3) = x^3$, which follow from~\cite{KaSo2}, Paragraph~2.5. The form of~$\sbb(v)$ can be derived in the same way, but the computation is more tedious: Using the formula for the coaction of~$y = \eta_1$ given in Paragraph~\ref{YetDrinfFirst} and the equations for~$\sa(\eta_i)$ from~\cite{KaSo2}, loc.~cit. yields
\begin{align*}
\sbb \left( v \right)
= {} & \frac{1}{8} (h_1+h_2+h_3+h_4)  \left(d_1 + \zeta^3 d_2 + d_3 - \zeta^3 d_4 \right) \\
&+ \frac{1}{8}(h_1+\zeta^2 h_2-h_3-\zeta^2 h_4) \left(\zeta^3d_1 + d_2 - \zeta^3 d_3 + d_4 \right)\\
&+\frac{1}{8} (h_1-h_2+h_3-h_4)  \left(d_1 - \zeta^3 d_2 + d_3 + \zeta^3 d_4 \right)  \\
&+ \frac{1}{8}(h_1-\zeta^2 h_2-h_3+\zeta^2 h_4)  \left(- \zeta^3d_1 +  d_2 + \zeta^3 d_3 + d_4 \right)
\end{align*}
The assertion now follows by collecting like terms, writing the result in terms of the generators, and using their commutation relations.
\qed
\end{pf}
\enlargethispage{10pt}

\subsection[The center]{} \label{Cent}
In order to determine the Wedderburn decomposition of~$B$, it will be helpful to know the center of~$B$. It is described in the following theorem:
\pagebreak
\begin{samepage}
\begin{thm}
The center~$Z(B)$ of~$B$ has dimension~$11$. The four sets of vectors
\begin{enumerate}
\item
$c_1, \quad c_4, \quad c_2 + c_3, \quad d_1 + d_2 + d_3 + d_4$

\item
$(c_1 + c_2 + c_3 + c_4) h_2, \quad (d_1 + d_2 + d_3 + d_4) h_2$

\item
$(c_1 + c_4) h_3, \quad (c_2 + c_3) h_3, \quad (d_1 + d_2 + d_3 + d_4) h_3$

\item
$(c_1 + c_2 + c_3 + c_4) h_4, \quad (d_1 + d_2 + d_3 + d_4) h_4$
\end{enumerate}
together form a basis of~$Z(B)$.
\end{thm}
\end{samepage}
\begin{pf}
\begin{pflist}
\item
Suppose that $b = \sum\displaystyle_{j=1}^{4} a_{j} \star g_j$ lies in the center of~$B$. That~$b$~commutes with elements of the form $h = \A \star g$ for some $g \in G$ means that
\[\sum_{j=1}^{4}a_{j}\star g_j g= bh=hb =\left( \A \star g\right)\sum_{j=1}^{4}a_{j}\star g_j = \sum_{j=1}^{4}(g .a_{j})\star g_j g\]
Thus, for each $j = 1,2,3,4$, we must have $g.a_{j} = a_j$, and therefore $a_j$ lies in the space of invariants, which we denote by~$A^H$.

\item
On the other hand, $b$~commutes with elements of the form~$b' = a' \star \HH$ provided that
\[\sum_{j=1}^{4} a' a_{j} \star g_j = b' b = b b' = 
\left(\sum_{j=1}^{4} a_{j} \star g_j \right) \left( a' \star \HH \right) = 
\sum_{j=1}^{4} (g_j.a') a_{j} \star g_j\]
This condition is equivalent to the condition
\begin{equation*}
a' a_j = (g_j.a') a_j
\end{equation*}
for all $j = 1,2,3,4$.

\item
In fact, it is sufficient to test this condition in the two cases~$a'=x$ and \mbox{$a'=y$}, because~$b$ will be contained in the center if and only if it commutes with the four generators~$u$, $v$, $r$, and~$s$ introduced at the beginning of Paragraph~\ref{AlgBiprod}. The generators~$r$ and~$s$ are of the form treated in Step~(1), and the generators~$u$ and~$v$ correspond to the two cases~$a'=x$ and~$a'=y$.

\item
Discussing each~$j$ separately, we see that the condition $a' a_{j} = (g_j .a') a_{j}$ is vacuous for~$j=1$. For~$j=2$, the condition states that
\[x a_2 = (g_2.x) a_2 = x^3 a_2 \qquad \text{and} \qquad y a_2 = (g_2.y) a_2 = x^3 y a_2\]
in other words, $a_2 = x^2 a_2$ and $a_2 = x^3 a_2$, which means that~$x a_2 = a_2$.

\begin{samepage}
For~$j=3$, the condition states that
\[x a_3 = (g_3.x) a_3 = x a_3 \qquad \text{and} \qquad y a_3 = (g_3.y) a_3 = x^2 y a_3\]
The first part is vacuous and the second is equivalent to the condition $x^2 a_3 = a_3$.

\end{samepage}

For~$j=4$, the condition states that
\[x a_4 = (g_4.x) a_4 = x^3 a_4 \qquad \text{and} \qquad y a_4 = (g_4.y) a_4 = x y a_4\]
in other words, $a_4 = x^2 a_4$ and $a_4 = x a_4$. Obviously, the second part of the condition implies the first.

\item
The group-like elements of~$A$ are divided into four orbits under the action of~$G$, and~$A^H$ is spanned by the sum of the elements in each of these orbits, i.e., is spanned by the four elements
\[\omega_1 \qquad \qquad \omega_4 \qquad \qquad \omega_2 + \omega_3 \qquad  \qquad
\eta_1 + \eta_2 + \eta_3 + \eta_4  \]
It is immediate from the definition of these elements that $\omega_1 = \A$,
\mbox{$\omega_4 = x^2$}, \mbox{$\omega_2 + \omega_3 = x + x^3$}, and
$\eta_1 + \eta_2 + \eta_3 + \eta_4 = (\A + x + x^2 + x^3) y$. We can therefore expand~$a_j$ in the form
\[a_{j} = \lambda_{1,j} \A + \lambda_{2,j} x^2 + \lambda_{3,j} (x + x^3)
+ \lambda_{4,j} (\A + x + x^2 + x^3) y\]

For $j=1$, we saw above that there are no extra conditions on~$a_1$, so the coefficients~$\lambda_{k,1}$ can be arbitrary.
For~$j=2$, the condition~$x a_2 = a_2$ implies that
$\lambda_{1,2} = \lambda_{2,2} = \lambda_{3,2}$, and so~$a_2$ is a linear combination of the vectors
\[\A + x + x^2 + x^3 \qquad \text{and} \qquad (\A + x + x^2 +x^3) y\]
For~$j=3$, the condition~$x^2 a_3 = a_3$ implies only that
$\lambda_{1,3} = \lambda_{2,3}$, and so~$a_3$ is a linear combination of the three vectors
\[\A + x^2 \qquad \qquad x + x^3 \qquad \qquad (\A + x + x^2 +x^3) y\]
The condition on~$a_4$ is the same as the condition on~$a_2$, and so~$a_4$ is a linear combination of the same vectors as~$a_2$, which can, as we have already seen, also be written as
\[\omega_1 + \omega_2 + \omega_3 + \omega_4 \qquad \text{and} \qquad
\eta_1 + \eta_2 + \eta_3 + \eta_4\]
Rewriting these results in terms of~$c_i$ and~$d_i$ instead of~$\omega_i$ and~$\eta_i$ yields the assertion, where the four items in the list correspond to the index values~$j=1$, $j=2$, $j=3$, and~$j=4$ in the preceding discussion.
\qed
\end{pflist}
\end{pf}

\section{Representations} \label{Sec:Rep}
\subsection[Group-like elements]{} \label{Group}
As we pointed out in Paragraph~\ref{Biprod}, the dual space~$A^*$ is a Yetter-Drinfel'd Hopf algebra over~$H^*$, and the corresponding biproduct~$A^*\star H^*$ is isomorphic to~$B^*$. The Hopf
algebra~$H^* = K[G]^*$ is itself a group algebra, namely the group algebra of the character group~$\hat{G}$, whose elements are considered as elements of~$H^*$ by linear extension. In our situation, the nondegenerate symmetric bicharacter~$\theta$ described in Paragraph~\ref{YetDrinfFirst} yields a specific bijection between~$G$ and~$\hat{G}$: For~$i=1, 2, 3, 4$, we define the character~$\gamma_i$ by
\[\gamma_i(g_j) = \theta (g_i, g_j)\]
for all~$j=1, 2, 3, 4$. We then have
$\hat{G} = \{\gamma_1, \gamma_2, \gamma_3, \gamma_4\}$, where~$\gamma_1 = \eh$.

The action of~$\gamma_i$ on~$A^*$ is closely related to the action of~$g_i$ on~$A$:
For~$\varphi \in A^*$ and~$a \in A$, we have by Equation~(\ref{dualaction}) and \cite{KaSo2}, Lemma~4.3 that
\begin{equation} \label{dualaction2}
\left( \gamma_i .\varphi \right) (a) =  \left( \gamma_i \ot \varphi \right)\delta(a)
= \gamma_i(a^\1) \varphi (a^\2) = \varphi(g_i^{-1}.a)
\end{equation}
where of course~$g_i = g_i^{-1}$ in our situation.

This equation will be used below to determine which of the group-like elements in~$B^*$ are central. But first, we record the group-like elements of~$B$:
\begin{lem} \label{GrouplB}
We have
\begin{align*}
{\bf G}(B) &=  \left\langle c_4 \right\rangle \times \left\langle h_2 \right\rangle \times \left\langle h_3 \right\rangle \cong \Z_2 \times \Z_2 \times \Z_2
\end{align*}
Moreover, $c_4$ is the only nontrivial central group-like element of~$B$.
\end{lem}
\begin{pf}
As we already mentioned in Paragraph~\ref{CoalgBiprod}, the group-like elements of~$B$ are of the form $a \star g_j$, where $a \in A$ is group-like and coinvariant. Of the eight group-like elements
$\omega_1, \omega_2, \omega_3, \omega_4, \eta_1, \eta_2, \eta_3, \eta_4$ of~$A$, only~$\omega_1 = \A$ and~$\omega_4$ are coinvariant. Because these are also invariant, all the group-like elements of~$B$ commute, so that they form the elementary abelian group of order~$8$ stated above. As we have seen explicitly in the proof of Theorem~\ref{Cent} that~$h_2$, $h_3$, and~$h_4$ are not central, $c_4$~is the only nontrivial central group-like element of~$B$.
\qed
\end{pf}

Alternatively, this lemma can be stated in terms of the generators introduced at the beginning of Paragraph~\ref{AlgBiprod}, as we have $c_4 = u^2$, $h_2 = r$, and~$h_3 = s$.

The group of group-like elements of~$B^*$ has a very similar form:
\begin{prop} \label{GrouplB*}
There are three elements $\chi_1, \chi_2, \chi_3 \in {\bf G}(B^*)$ that satisfy
\begin{enumerate}
\item
$\chi_1(d_i) = -1$ and $\chi_1(c_i) = \chi_1(h_i) = 1$

\item
$\chi_2 (h_2) = \chi_2 (h_4)=-1 $ and $\chi_2 (c_i) = \chi_2(d_i) = \chi_2(h_3) =1$

\item
$\chi_3(h_3) = \chi_3(h_4) = -1$ and $\chi_3(c_i) = \chi_3(d_i) = \chi_3(h_2) =1$
\end{enumerate}
for all $i = 1, 2, 3, 4$. With these elements, we have
\begin{align*}
{\bf G}(B^*) &= \left\langle \chi_1 \right\rangle \times\left\langle\chi_2 \right\rangle \times \left\langle \chi_3 \right\rangle\cong \Z_2 \times \Z_2\times \Z_2
\end{align*}
Moreover, $\chi_1$ is the only nontrivial central group-like element of $B^*$.
\end{prop}
\begin{pf}
\begin{pflist}
\item
Since  $B^* \cong A^*\star H^*$, it follows as in the preceding lemma that the group-like elements of~$B^*$ are of the form $\rho \star \gamma$, where $\gamma \in \hat G$ and $\rho \in {\bf G}(A^*)$ is $H^*$-coinvariant. In view of Equation~(\ref{dualcoaction}), a linear form is  $H^*$-coinvariant if and only if it is an $H$-linear map to the base field, and similarly Equation~(\ref{dualaction}) implies that it is $H^*$-invariant if and only if it is an $H$-colinear map to the base field. On the other hand, Equation~(\ref{dualaction2}) above shows that a linear form is $H^*$-invariant if and only if it is $H$-linear, so that all these four properties are equivalent.

\item
According to \cite{KaSo2}, Paragraph~2.1, there are eight one-dimensional characters of~$A$, which were denoted there by $\varepsilon_1$, $\varepsilon_2$, $\varepsilon_3$, $\varepsilon_4$,
$\rho_1$, $\rho_2$, $\rho_3$, and~$\rho_4$. Their definition depended on the exact order of~$\zeta$, but the only $H$-linear ones were~$\el_1$ and~$\el_4$, where~\mbox{$\el_1 = \ea$} was the counit of~$A$, and~$\el_4$ was in each of the cases determined by the condition
that~$\el_4(x) = \K$ and~$\el_4(y) = -\K$. Thus~$B^*$ contains eight group-like elements, namely
\[\el_1 \star \gamma_1 \mspace{28mu} \el_1 \star \gamma_2 \mspace{28mu} \el_1 \star \gamma_3 \mspace{28mu} \el_1 \star \gamma_4 \mspace{28mu}
\el_4 \star \gamma_1 \mspace{28mu} \el_4 \star \gamma_2 \mspace{28mu} \el_4 \star \gamma_3 \mspace{28mu} \el_4 \star \gamma_4\]
The $H^*$-invariance of~$\el_1$ and~$\el_4$ and the multiplication formula in a smash product stated in Equation~(\ref{Multip}) together imply that these elements commute; since
\mbox{$\el_4^2 = \ea$}, they all have at most order~$2$.

\item
We set $\chi_1 \deq \el_4 \star \gamma_1$ and then have $\chi_1(d_i) = -1$
and~$ \chi_1(c_i) = \chi_1(h_i) = 1$ as required. By letting $\chi_2 \deq \el_1 \star \gamma_3$, the requirements \mbox{$\chi_2(h_2) = \chi_2(h_4) = -1$}
and~$\chi_2(c_i) = \chi_2 (d_i) = \chi_2(h_3) = 1$ are also satisfied. But for~$\gamma_2$, we have according to the definition that~$\gamma_2(g_2) = \zeta^2$ and~$\gamma_2(g_3) = -1$, while~\mbox{$\gamma_4(g_2) = -\zeta^2$} and~\mbox{$\gamma_4(g_3) = -1$}. So if~$\zeta^2 = 1$, we set $\chi_3 \deq \el_1 \star \gamma_2$, while we set~$\chi_3 \deq \el_1 \star \gamma_4$ if ~$\zeta^2 = -1$. It is obvious that these three characters generate~${\bf G}(B^*)$.

\item
Because~$\el_4$ is $H^*$-invariant, we have for $\varphi \in A^*$ and~$\gamma \in H^*$ that
\[(\varphi \star \gamma) \chi_1 = \varphi (\gamma_\1.\el_4) \star \gamma_\2
= \varphi \el_4 \star \gamma = \chi_1 (\varphi \star \gamma) \]
so that~$\chi_1$ is central. On the other hand, we have for  $i= 2, 3, 4$ that
\[(\ea \star \gamma_i) (\varphi \star \eh) =  (\gamma_i.\varphi) \star \gamma_i
\qquad \text{but} \qquad (\varphi \star \eh) (\ea \star \gamma_i) = \varphi \star \gamma_i \]
Since in general~$\gamma_i.\varphi \neq \varphi$ by Equation~(\ref{dualaction2}) above,
$\ea \star \gamma_i$ is not central. This shows that~$\chi_1$ is the only nontrivial central group-like element in~$B^*$.
\qed
\end{pflist}
\end{pf}

At least part of the preceding proposition can also be derived with the help of the presentation of~$B$ given in Proposition~\ref{AlgBiprod}. The values of our characters on the generators~$u$, $v$, $r$, and~$s$ introduced there are given in the following table:
\begin{center}
\renewcommand\arraystretch{1.5}
\begin{tabular}{|c|c|c|c|c|} \hline
& $u$ & $v$ & $r$ & $s$ \\ \hline
$\chi_1$ & $1$ & $-1$ & $1$ & $1$ \\ \hline
$\chi_2$ & $1$ & $1$ & $-1$ & $1$ \\ \hline
$\chi_3$ & $1$ & $1$ & $1$ & $-1$ \\ \hline
\end{tabular}
\end{center}

We also note that the existence of nontrivial central group-like elements in~$B$ and~$B^*$ is consistent with~\cite{MaPn}, Theorem~1, which states under our assumptions on the base field that a semisimple Hopf algebra of prime power dimension always contains such elements.

\subsection[The Wedderburn decomposition]{} \label{Wedderb}
We have already mentioned in Paragraph~\ref{Biprod} that biproducts are semisimple and cosemisimple if their two factors are. Here, this implies that the biproduct~$B$ under consideration is semisimple and cosemisimple, because the Yetter-Drinfel'd Hopf algebra~$A$ is semisimple by~\cite{KaSo2}, Proposition~2.1 and cosemisimple by~\cite{KaSo2}, Proposition~2.4. But we can determine the Wedderburn decomposition of~$B$ and~$B^*$ precisely:

\begin{thm} \label{RepsBB*}
Both $B$ and~$B^*$ have~eight $1$-dimensional,~two $2$-dimensional, and one~$4$-dimensional irreducible representations.
\end{thm}
\begin{pf}
Since~$\dm Z(B)=11$ by Theorem~\ref{Cent} and~$\left|{\bf G}(B^*)\right|=8$ by Proposition~\ref{Group}, there are exactly three irreducible representations of~$B$ of dimension larger than~$1$. The sum of the squares of the dimensions of these three representations is~$32-8=24$, which implies that~$B$ has, up to isomorphism, two $2$-dimensional and one $4$-dimensional irreducible representations.

To determine the Wedderburn decomposition of~$B^*$, recall from Paragraph~\ref{Biprod} that~$B^*$ is the smash product of~$A^*$ and~$H^* = K[\hat G]$ as an algebra. Therefore, the correspondence between the irreducible representations of~$B^*$ and the irreducible representations of~$A^*$ is described by Clifford theory, which in our situation yields the following form of the irreducible representations: Because the coalgebra~$A$ has a basis consisting of group-like elements, $A^*$~is a commutative semisimple algebra whose simple modules are one-dimensional and given by evaluation at one of these group-like elements. Its primitive idempotents are therefore the dual basis elements
\[p_1 \qquad p_2 \qquad p_3 \qquad p_4 \qquad q_1 \qquad q_2 \qquad q_3 \qquad q_4\]
of the basis~$\omega_1, \omega_2, \omega_3, \omega_4, \eta_1, \eta_2, \eta_3, \eta_4$. The action of~$\hat G$ on~$A^*$ permutes these primitive idempotents. According to the theory, we have to choose a representative~$p$ in each orbit, consider the stabilizer~$T_p \subset \hat G$, often called the inertia group, of this representative, and find the simple \mbox{$T_p$-modules}. For each such simple \mbox{$T_p$-module~$V$}, the induced module~$\Ind^{\hat G}_{T_p} V$ can be endowed with a \mbox{$B^*$-module} structure that makes it simple (cf.~\cite{CR2}, Proposition~(11.16); \cite{KMM},~Corollary~3.5). This yields a bijective correspondence between the isomorphism classes of simple \mbox{$T_p$-modules} in the various orbits and the isomorphism classes of simple \mbox{$B^*$-modules}, the so-called Clifford correspondence.

It follows from Equation~(\ref{dualaction2}) that~$\hat G$ permutes the primitive idempotents in~$A^*$ in the same way as~$G$ permutes the group-like elements in~$A$ (cf.~\cite{KMM}, Equation~(2.3)). From Paragraph~\ref{YetDrinfFirst}, we know that the orbits of the \mbox{$G$-action} on the set of group-like elements of~$A$ are $\{\omega_1\}$, $\{\omega_4\}$, $\{\omega_2, \omega_3\}$, and~$\{\eta_1, \eta_2, \eta_3, \eta_4\}$. Therefore, the orbits of the \mbox{$\hat G$-action} on the set of primitive idempotents in~$A^*$ are $\{p_1\}$, $\{p_4\}$, $\{p_2, p_3\}$, and~$\{q_1, q_2, q_3, q_4\}$.

Since~$T_{p_1} = T_{p_4} = \hat G$, the Clifford correspondence yields four $1$-dimensional irreducible representations for each of the orbits~$\{p_1\}$ and~$\{p_4\}$, so eight \mbox{$1$-dimen}\-sional irreducible representations in total. Since~$T_{q_1}= \{1_{\hat G}\}$, the Clifford correspondence yields one $4$-dimensional irreducible representation for the
orbit~$\{q_1, q_2, q_3, q_4\}$. Finally, since~$g_3$ fixes~$\omega_2$, it follows from Equation~(\ref{dualaction2}) that~$\gamma_3$ fixes~$p_2$, so
that~$T_{p_2} = \{1_{\hat G}, \gamma_3\}$. The two $1$-dimensional representations of~$T_{p_2}$ therefore give rise to two $2$-dimensional irreducible representations of~$B^*$.
\qed
\end{pf}

In Paragraph~\ref{GrothB}, we will describe in addition the product in the Grothendieck ring of~$B$.

\section{Hopf subalgebras of dimension~$16$} \label{Sec:2by16}
\subsection[Constructing the Hopf subalgebra]{} \label{N16}
The algebra~$C$ introduced in Paragraph~\ref{YetDrinfFirst} is a Yetter-Drinfel'd Hopf subalgebra of~$A$, and therefore the biproduct~$B = A \star H$ contains the biproduct
\[N \deq C \star H\]
as a Hopf subalgebra. Clearly, $N$ has dimension~$16$, and the elements~$c_i h_j$ for~$i,j = 1,2,3,4$ form a basis of~$N$. This implies that~$N$ is generated as an algebra by~$c_2$, $c_4$, $h_2$, and~$h_3$. This generating set has the advantage of being close to the group-like elements of~$A$. Alternatively, $N$ is generated by the elements~$u$, $r$, and~$s$ introduced in Paragraph~\ref{AlgBiprod}, which satisfy the relations
\[u^4 = 1 \qquad r u = u^3 r \qquad s u = u s \qquad r^2 = 1 \qquad rs = sr \qquad s^2 = 1 \]
from Proposition~\ref{AlgBiprod}. An argument that is very similar to the one used there shows that these relations are defining. Comparing this presentation to the usual presentation of the dihedral group (cf.~\cite{SerRep}, Section~5.3), we see that~$N$ is isomorphic to the
group ring~$K[D_8 \times \Z_2]$ as an algebra, but not as a coalgebra, as it is not cocommutative, as for example the formula for the coproduct of~$c_2$ in Paragraph~\ref{CoalgBiprod} shows. This issue will be revisited at the end of Paragraph~\ref{QN16}.

The groups of group-like elements of~$N$ and~$N^*$ have the following form:
\begin{lem}
We have
\begin{align*}
{\bf G}(N) &=  \left\langle c_4 \right\rangle \times \left\langle h_2 \right\rangle \times \left\langle h_3 \right\rangle\cong \Z_2 \times \Z_2 \times \Z_2
\end{align*}
and also ${\bf G}(N^*) \cong \Z_2 \times \Z_2 \times \Z_2$.
\end{lem}
\begin{pf}
Since all the group-like elements of~$B$ are already contained in~$N$, we have
${\bf G}(N) = {\bf G}(B)$, so that the first assertion follows from Lemma~\ref{Group}.

As in the proof of Proposition~\ref{Group}, we have $N^* \cong C^* \star H^*$, and under this isomorphism a group-like element~$\chi \in {\bf G}(N^*)$ decomposes in the form
$\rho \star \gamma$, where $\gamma \in \hat{G}$ and $\rho \in {\bf G}(C^*)$ is $H^*$-coinvariant. We also discussed there that $H^*$-coinvariance is equivalent to $H^*$-invariance as well as to $H$-linearity and \mbox{$H$-colinearity}. If~$\chi' \in {\bf G}(N^*)$ is another group-like element that decomposes in the form~$\rho' \star \gamma'$, the $H^*$-invariance implies that~$\chi \chi'$ corresponds to~$\rho \rho' \star \gamma \gamma'$.

There are clearly four algebra homomorphisms from~$C$ to~$K$, determined by the image of~$x$, which must be a fourth root of unity. Out of these, two are $H$-linear,
namely~$\varepsilon'_1 \deq \varepsilon_C$, the counit, which satisfies
$\varepsilon'_1(x) = \K$, and the character~$\varepsilon'_2$ determined by
$\varepsilon'_2(x) = -\K$. It also satisfies 
$\varepsilon'_2(\omega_2) = \varepsilon'_2(\omega_3) = -\K$, so that 
\mbox{$\varepsilon_2'^2 = \varepsilon_C$}. This shows that $|{\bf G}(N^*)|=8$ and
\mbox{${\bf G}(N^*) \cong \Z_2 \times \Z_2 \times \Z_2$}.
\qed
\end{pf}

In analogy with Paragraph~\ref{Group}, we choose generators~$\chi'_1$, $\chi'_2$, and $\chi'_3$ of~${\bf G}(N^*)$ whose values on the generators of~$N$ are given by the following table:
\begin{center}
\renewcommand\arraystretch{1.5}
\begin{tabular}{|c|c|c|c|} \hline
& $u$ & $r$ & $s$ \\ \hline
$\chi'_1$ & $-1$ & $1$ & $1$ \\ \hline
$\chi'_2$ & $1$ & $-1$ & $1$ \\ \hline
$\chi'_3$ & $1$ & $1$ & $-1$ \\ \hline
\end{tabular}
\end{center}

Explicitly, $\chi'_1$ corresponds to $\varepsilon'_2 \star \gamma_1$, while~$\chi'_2$ corresponds to $\varepsilon'_1 \star \gamma_3$, and~$\chi'_3$~corresponds to $\el'_1 \star \gamma_2$ if~$\zeta^2 = 1$ and to~$\el'_1 \star \gamma_4$ if~$\zeta^2 = -1$.

The preceding lemma is consistent with a general fact about semisimple Hopf algebras of dimension~$16$ that are neither commutative nor cocommutative: Such a Hopf algebra has an abelian group of group-likes of order~$8$ if and only if the dual Hopf algebra does (cf.~\cite{Ka}, Proposition~3.1).

The isomorphism ${\bf G}(N^*) \cong \Z_2 \times \Z_2 \times \Z_2$ can also be established in a more direct way that does not depend on the biproduct perspective:
\begin{cor}
If $\Z_2 = \{1,-1\}$ is written multiplicatively, the map
\[{\bf G}(N^*) \to \Z_2 \times \Z_2 \times \Z_2,~\chi \mapsto (\chi(u), \chi(r), \chi(s))\]
is an isomorphism of groups.
\end{cor}
\begin{pf}
For~$\chi \in {\bf G}(N^*)$, the relation $\chi(r u) = \chi(u^3 r)$ implies
that~$\chi(u)^2 = 1$. This condition and the relations
$r^2 = s^2 = 1$ imply that~$\chi$ must map all three generators~$u$, $r$, and~$s$ to~$\pm 1$.
So the above map is well-defined, and also injective, because~$\chi$ is determined by its values on the generators. On the other hand, it is not complicated to check that the above defining relations of~$N$ are satisfied by an arbitrary choice of these signs, so the map is also surjective.  Our assertion will therefore follow if we can show that it is a group homomorphism.

For~$\chi, \chi' \in {\bf G}(N^*)$, we clearly have $(\chi \chi')(r) = \chi(r) \chi'(r)$, since~$r$ is group-like, and the same equation holds for~$s$. For~$u$, we have by Proposition~\ref{CoalgBiprod} that
\begin{align*}
(\chi \chi')(u) &= \frac{\K}{2} \left(\chi(u) \chi'(u) + \chi(u) \chi'(u)^3
+ \chi(u)^3 \chi(s) \chi'(u) - \chi(u)^3 \chi(s) \chi'(u)^3 \right) \\
&= \chi(u) \chi'(u)
\end{align*}
so that the requirement is also fulfilled for this generator.
\qed
\end{pf}

We note that the fact that~$\chi(u) = \pm 1$ recorded in the preceding proof implies that
\[\chi(u) = \chi(c_2) = \chi(c_3)\]
for all~$\chi \in {\bf G}(N^*)$.

Of the eight group-like elements of~$N^*$, four are central:
\begin{prop}
The group~${\bf G}(N^*) \cap Z(N^*)$ is isomorphic to $\Z_2 \times \Z_2$ and is generated by
$\chi'_1$ and~$\chi'_2$.
\end{prop}
\begin{pf}
If~$\chi \in {\bf G}(N^*)$ corresponds to~$\rho \star \gamma \in C^* \star H^*$ as in the proof of the preceding lemma, we have for~$\varphi \in N^*$ in view of the $H$-colinearity of~$\rho$ that
\begin{align*}
(\varphi \chi) (a \star h) &=
\varphi (a_\1 \star a_\2{}^\1 h_\1) \; \rho(a_\2{}^\2) \gamma(h_\2) \\
&= \varphi (a_\1 \star h_\1) \; \rho(a_\2) \gamma(h_\2)
\end{align*}
on the one hand and
\begin{align*}
(\chi \varphi) (a \star  h) =
\rho(a_\1) \gamma(a_\2{}^\1 h_\1) \; \varphi(a_\2{}^\2 \star h_\2)
\end{align*}
on the other hand, so that~$\chi$ is central if and only if
\[a_\1 \star h_\1 \; \rho(a_\2) \gamma(h_\2) =
\rho(a_\1) \gamma(a_\2{}^\1 h_\1) \; a_\2{}^\2 \star h_\2 \]
for all $a \in C$ and all~$h \in H$.

If~$\chi = \chi'_1$, we have $\rho = \varepsilon'_2$ and $\gamma = \gamma_1 = \eh$, so that this condition becomes
\[a_\1 \star h \; \varepsilon'_2(a_\2) = \varepsilon'_2(a_\1)  \; a_\2 \star h \]
which is satisfied since~$C$ is cocommutative. So $\chi'_1$ is central.

If~$\chi = \chi'_2$, we have $\rho = \varepsilon'_1 = \varepsilon_C$ and $\gamma = \gamma_3$,  so that this condition becomes
\[a \star h_\1 \gamma_3(h_\2) =  \gamma_3(a^\1 h_\1) \; a^\2 \star h_\2 \]
which will hold if and only if~$\gamma_3(a^\1) a^\2 = a$ for all $a \in C$. But because
$\gamma_3(g_3) = \theta(g_3, g_3) = 1$, the formulas for~$\delta(x)$ and~$\delta(x^3)$ given in Paragraph~\ref{YetDrinfFirst} show that this is indeed the case. So $\chi'_2$ is also central.

If~$\chi = \chi'_3$, we have $\rho = \varepsilon'_1 = \varepsilon_C$ and $\gamma = \gamma_2$ or $\gamma = \gamma_4$, depending on the order of~$\zeta$. But if~$a = \omega_2$, the formula for~$\delta(\omega_2)$ given in Paragraph~\ref{YetDrinfFirst} yields that $\gamma_2(a^\1) a^\2 = \gamma_4(a^\1) a^\2 =  \omega_3$. So the condition is not met in both cases and~$\chi'_3$ is not central.
\qed
\end{pf}

\subsection[Subquotients]{} \label{QN16}
A characteristic property of~$N$ is the following:
\begin{prop} \label{3quotients}
The Hopf algebra~$N$ has exactly three quotient Hopf algebras of dimension~$8$. All of these quotients are cocommutative and isomorphic to either $K[\Z_2 \times \Z_2 \times \Z_2]$ or~$K[D_8]$.
\end{prop}
\begin{pf}
\begin{pflist}
\item
If~$N$ has a quotient Hopf algebra~$F$ of dimension~$8$, then~$N^*$ has a Hopf subalgebra~$F^*$ that has index~$2$ and is therefore necessarily normal by~\cite{KoMa}, Proposition~2. Thus~$N^*$ fits into an extension of the form
\[F^* \hookrightarrow N^* \twoheadrightarrow M\]
and so~$N$ fits into an extension of the form
\[M^* \hookrightarrow N \twoheadrightarrow F\]
where~$M^*$ is normal and has dimension~$2$. As pointed out in Paragraph~\ref{Extens}, $M^*$~is then central, and therefore~$M^*=K\langle g \rangle$ for a central group-like element~$g$ of order~$2$, so that~$F \cong N/( (M^*)^+ N)$.

\item
It follows from Lemma~\ref{N16} that~$N$ contains exactly three central group-like elements of order~$2$, namely~$c_4$, $h_3$, and~$c_4 h_3$. Therefore it has exactly three quotient Hopf algebras of dimension~$8$, namely
\[F_1 \deq N/\left( K\langle c_4 \rangle^+ N\right) \qquad
F_2 \deq N/\left( K \langle h_3 \rangle^+ N \right) \qquad
F_3 \deq N/\left( K \langle c_4 h_3 \rangle^+ N \right) \]
In all three of these quotient spaces, the images of the group-like elements~$r=h_2$ and~$s=h_3$ are group-like.

\item
In~$F_1$, we have $\bar{c}_4 = \bar{u}^2 = 1_{F_1}$, so that~$F_1$ is generated by $\bar{c}_2$, $\bar{h}_2$, and~$\bar{h}_3$, or alternatively by~$\bar{u}$, $\bar{r}$, and~$\bar{s}$, since 
$\bar{u} = \bar{c}_2$ in~$F_1$.
We have seen in Proposition~\ref{CoalgBiprod} that
\begin{align*}
\db(u) &= \frac{\K}{2} (u \ot u + u \ot u^3 +  u^3 s \ot u - u^3 s \ot u^3)
\end{align*}
In~$F_1$, this equation reduces to~$\Delta_{F_1}(\bar{u}) = \bar{u} \ot \bar{u}$. In fact, the relations recorded at the beginning of Paragraph~\ref{N16} imply that \mbox{$\bar{u}^2 = \bar{r}^2 = \bar{s}^2 = 1$} and that these three generators commute in~$F_1$, which shows that~$F_1 \cong K[\Z_2 \times \Z_2 \times \Z_2]$.

\item
In~$F_2$, we have $\bar{h}_3 = \bar{s} = 1_{F_2}$, so that~$F_2$ is generated by $\bar{c}_2$, $\bar{c}_4$, and~$\bar{h}_2$, or alternatively by~$\bar{u}$ and~$\bar{r}$.
We have seen in Paragraph~\ref{CoalgBiprod} that
\begin{align*}
\db(c_2) &= \frac{1}{2} c_2(\HH + h_3) \ot c_2 + \frac{1}{2} c_2 (\HH - h_3) \ot c_3
\end{align*}
In~$F_2$, this reduces to~$\Delta_{F_2}(\bar{c}_2) = \bar{c}_2 \ot \bar{c}_2$. Now we know from Equation~(\ref{c2h2}) that $(c_2 h_2)^2 = c_4$, so that $(c_2 h_2)^3 = c_3 h_2$ and $(c_2 h_2)^4 = 1$. This shows that~$F_2$ is already generated by $\bar{c}_2$ and~$\bar{h}_2$ alone. Also, we have
\[(c_2 h_2)^3 h_2 = c_3 = h_2 (c_2 h_2)\]
In view of the presentation of the dihedral group~$D_8$ already mentioned in Paragraph~\ref{N16}, this implies that~$F_2 \cong K[D_8]$.

\item
In~$F_3$, we have~\mbox{$\bar{u}^2 = \bar{s}$}, so that~$F_3$ is generated by $\bar{u}$ and~$\bar{r}$. The formula for~$\db(u)$ already recalled above gives~$\Delta_{F_3}(\bar{u}) = \bar{u} \ot \bar{u}$ also in this case. The relations recorded at the beginning of Paragraph~\ref{N16} yield directly that~\mbox{$\bar{u}^4 = 1$}
and~$\bar{r} \bar{u} = \bar{u}^3 \bar{r}$, so that also~$F_3 \cong K[D_8]$. As in the previous step, $F_3$ is also already generated by $\bar{c}_2$ and~$\bar{h}_2$ alone.
\qed
\end{pflist}
\end{pf}

From the classification results in~\cite{Ka},~Table 1, we know that there are exactly two Hopf algebras of dimension~$16$ for which both the group of group-likes and the group of group-likes of the dual are elementary abelian of order~$8$, namely those denoted there by~$H_{d:1,1}$
and~$H_{d:-1, 1}$. The Hopf algebra~$H_{d:-1, 1}$ is isomorphic
to~\mbox{$H_8 \ot K\Z_2$}, where~$H_8$ is the Kac-Paljutkin Hopf algebra of dimension~$8$. Since~$H_8 \ot K\Z_2$ has a noncocommutative quotient isomorphic to~$H_8$, the proposition above implies that~$N$ must be isomorphic to~$H_{d:1,1}$. By \cite{Ka}, Remark~1.4, we then know that~$N$ is isomorphic to~$K [D_8 \times \Z_2]_J$, which is triangular, and that~$N$ is self-dual. The twisting cocycle~$J$ is described in \cite{Ka}, Section~7, Example~2. In particular, we see again that~$N$ is isomorphic to~$K [D_8 \times \Z_2]$ as an algebra.

\subsection[Explicit isomorphisms]{} \label{H16}
In fact, we can construct an explicit isomorphism between~$N$ and~$H_{d:1,1}$. As described in~\cite{Ka}, Paragraph~3.2, Page~629ff, the Hopf algebra~$H_{d:1,1}$ is a bicrossed product~$K[\Gamma] \#^{\theta} K[L]$ for
$\Gamma =\left\langle x \right \rangle \times \left \langle y \right\rangle \times \left\langle z \right\rangle \cong \Z_2 \times \Z_2 \times \Z_2$
and \mbox{$L = \left\langle t \right\rangle \cong \Z_2$} in which the cocycle and the coaction are trivial. Here, and in the present paragraph only, we use the notation from~\cite{Ka}, loc.~cit., so that~$\theta$ does not denote the nondegenerate bicharacter from Paragraph~\ref{YetDrinfFirst}, but rather a dual cocycle like the one that was denoted by~$\kappa$ in Paragraph~\ref{Bicross}. Furthermore, $x$~and~$y$ do not denote the generators of~$A$ introduced in Paragraph~\ref{YetDrinfFirst}, but rather two of the three generators of~$\Gamma$ already appearing above. These generators give rise to the primitive idempotents~$e_{i,j,k}$ in~$K[\Gamma]$ via the formula
\[e_{i,j,k} \deq
\frac{1}{8} \left( 1 + (-1)^i x\right) \left( 1 + (-1)^j y\right) \left( 1 + (-1)^k z \right)\]
where the indices~$i$, $j$, and~$k$ each take the value~$0$ or~$1$. Conversely, the generators can be expressed in terms of the primitive idempotents via the formulas
\[x =\sum_{i,j,k=0}^{1} (-1)^i e_{i,j,k} \qquad
y =\sum_{i,j,k=0}^{1} (-1)^j e_{i,j,k} \qquad
z =\sum_{i,j,k=0}^{1} (-1)^k e_{i,j,k}\]
The action~$\rightharpoonup$ of the bicrossed product is determined by the facts that the unit element acts as the identity, while~$t$ acts as the algebra automorphism that takes the values
\[t \rightharpoonup x = y \qquad \qquad
t \rightharpoonup y = x \qquad \qquad
t \rightharpoonup  z = z \]
on the generators. Using the primitive idempotents, the dual cocycle~$\theta$ is defined as the map that sends~$t$ to
\[\theta(t) \deq \sum_{i,j,k,l,m,n=0}^{1} (-1)^{k (l+m)} e_{i,j,k} \ot e_{l,m,n} \]
and~$1$ to~$1 \ot 1$. Alternatively, we can express~$\theta(t)$ in terms of the generators. To do this, note that
\[ (-1)^{k (l+m)} = \frac{1}{2} (1 + (-1)^k + (-1)^{l+m} - (-1)^k (-1)^{l+m})\]
since both sides are equal to~$1$ when~$k=0$ and to~$(-1)^{l+m}$ when~$k=1$. Therefore we have
\begin{align*}
\theta (t) &= \frac{1}{2} \sum_{i,j,k,l,m,n}
\left(1 + (-1)^k + (-1)^{l+m} - (-1)^k (-1)^{l+m} \right) e_{i,j,k} \ot e_{l,m,n} \\
&=  \frac{1}{2}\left( \sum_{i,j,k} e_{i,j,k}+\sum_{i,j,k} (-1)^k e_{i,j,k}\right) \ot \sum_{l,m,n} e_{l,m,n} \\
&\quad +\frac{1}{2}\left(\sum_{i,j,k} e_{i,j,k} - \sum_{i,j,k} (-1)^k e_{i,j,k}\right) \ot \sum_{l,m,n} (-1)^{l+m} e_{l,m,n} \\
&=\frac{1}{2}\left( (1+z) \ot 1 +(1-z) \ot xy \right)
\end{align*}

The bicrossed product $H_{d:1,1} = K[\Gamma] \#^{\theta} K[L]$ just described is isomorphic to our Hopf subalgebra~$N$:
\begin{thm}
There exists a Hopf algebra isomorphism~$f \colon H_{d:1,1} \ra N$ satisfying
\begin{equation*} \label{isoHdN}
f(x) = h_2 \qquad \qquad
f(y) = c_4 h_2 \qquad \qquad
f(z) = h_3 \qquad \qquad
f(\bar t) = c_2
\end{equation*}
where~$a\in K[\Gamma ]$ is identified with~$a\# 1\in H_{d:1,1}$
and~$\bar{t} \deq 1 \# t \in H_{d:1,1}$.
\end{thm}
\begin{pf}
Because its cocycle is trivial by assumption, our crossed product is indeed a smash product, so that we can use the universal property of a smash product stated in Corollary~\ref{UnivPropCross}. Since both $\Gamma$ and~${\bf G}(N)$ are elementary abelian of order~$8$, there is an algebra homomorphism~$f_{K[\Gamma]} \colon K[\Gamma] \ra N$ defined by
\begin{equation*}
f_{K[\Gamma]}(x) = h_2 \qquad
f_{K[\Gamma]}(y) = c_4 h_2 \qquad
f_{K[\Gamma]}(z) = h_3
\end{equation*}
Similarly, as~$c_2$ has order~$2$, there is an algebra
homomorphism~\mbox{$f_{K[L]} \colon K[L] \ra N$} with~\mbox{$f_{K[L]} (t) = c_2$}. As we pointed out in Paragraph~\ref{UnivPropCross}, it is sufficient to verify the hypothesis of Corollary~\ref{UnivPropCross} on the generating set~$x,y,z$, where it holds since
\begin{align*}
f_{K[\Gamma]}(t \rightharpoonup x ) f_{K[L]}(t) & =f_{K[\Gamma]}(y)f_{K[L]}(t) = c_4 h_2 c_2 = c_2h_2 = f_{K[L]}(t) f_{K[\Gamma]}(x)\\
f_{K[\Gamma]}(t \rightharpoonup y )f_{K[L]}(t) &= f_{K[\Gamma]}(x)f_{K[L]}(t) = h_2 c_2 =
c_2 c_4 h_2 =  f_{K[L]}(t) f_{K[\Gamma]}(y)\\
f_{K[\Gamma]}(t \rightharpoonup z ) f_{K[L]}(t) &= f_{K[\Gamma]}(z) f_{K[L]}(t) = h_3 c_2 = c_2 h_3 = f_{K[L]}(t) f_{K[\Gamma]}(z)
\end{align*}
Thus Corollary~\ref{UnivPropCross} yields an algebra homomorphism~$f \colon H_{d:1,1} \to N$ that satisfies $f(a \# 1) = f_{K[\Gamma ]}(a)$ and~$f(1 \# c) = f_{K[L]}(c)$, and hence the asserted equations.

To see that~$f$ is a coalgebra homomorphism, it suffices to check the comultiplicativity of~$f$ on the generators of~$H_{d:1,1}$, namely~$x$, $y$, $z$, and~$\bar t$. Since~$x$, $y$, and~$z$ as well as their images are group-like, it is enough to verify this condition on~$\bar t$:
\begin{align*}
\left( f \ot f \right) \left(\Delta (\bar{t}) \right)
&= \left( f \ot f \right) \left(\frac{1}{2} \left( (1+z) \ot 1 + (1-z) \ot xy \right) \left(\bar{t} \ot \bar{t} \right) \right) \\
&=\frac{1}{2} \left( (\B + h_3) \ot \B  + (\B - h_3) \ot c_4 \right) \left(c_2 \ot c_2\right) \\
&= \db(c_2) = \db(f(\bar{t}))
\end{align*}
It is easily checked on the generators that~$f$ preserves the counit.

Since~$N$ is generated by~$h_2, c_4 h_2, h_3$, and~$c_2$, the map~$f$ is surjective. Therefore, since~$\dm H_{d:1,1} = \dm N = 16$, it is even bijective.
\qed
\end{pf}

The isomorphism~$f$ can be used to relate the one-dimensional representations of the two Hopf algebras. In particular, the one-dimensional characters~$\chi'_1$,~$\chi'_2$, and~$\chi'_3$ of~$N$ defined in Paragraph~\ref{N16} become the one-dimensional characters~\mbox{$\chi'_1 \circ f$}, $\chi'_2 \circ f$, and $\chi'_3 \circ f$ of~$H_{d:1,1}$, which were denoted in~\cite{Ka}, Page~633 by~$\chi$, $\varphi$, and~$\psi$, respectively.

\subsection[Uniqueness]{} \label{Proof2by16}
The sixteen-dimensional Hopf subalgebra~$N$ of~$B$ introduced in Paragraph~\ref{N16}, which is, as just discussed, isomorphic to $H_{d:1,1}$ and therefore a twisting of the group algebra of~$D_8 \times \Z_2$, is in fact the only such Hopf subalgebra:
\begin{thm}
The Hopf subalgebra~$N$ is the unique Hopf subalgebra of~$B$ of dimension~$16$. It is normal, but neither commutative nor cocommutative. The Hopf algebra~$B$ therefore fits into exactly one extension of the type
\[N \hookrightarrow B \twoheadrightarrow Z\]
with~$\dm N=16$ and~$\dm Z=2$. This extension is not abelian.
\end{thm}
\begin{pf}
Suppose that~$N'$ is a Hopf subalgebra of~$B$ of dimension~$16$. As in the proof of Proposition~\ref{QN16}, it follows from~\cite{KoMa}, Proposition~2 that~$N'$ is normal. Therefore, $B$~fits into the extension
\[N' \overset{\iota_{N'}}{\hookrightarrow} B \overset{\pi_Z}{\twoheadrightarrow} Z\]
where~$\iota_{N'}$ is the inclusion map and~$\pi_Z$ is the canonical projection to the
quotient~$Z \deq  B/B N'^+$ of dimension~$2$. By dualization, we obtain the extension
\[Z^* \overset{\pi_Z^*}{\hookrightarrow} B^* \overset{\iota_{N'}^*}{\twoheadrightarrow} N'^* \]
As pointed out in Paragraph~\ref{Extens}, the image~$\pi_Z^*(Z^*)$ is then central in~$B^*$, and is of course a group algebra by~\cite{M}, Theorem~2.3.1. But by Proposition~\ref{Group}, the only nontrivial central group-like element of~$B^*$ is~$\chi_1$, so that~$\pi_Z^*(Z^*)$ is spanned by~$\eb$ and~$\chi_1$. In particular, $\pi_Z^*(Z^*)$ together with its embedding into~$B^*$ is uniquely determined, so that dually~$B$ has a unique two-dimensional quotient Hopf algebra, namely~$Z$. But the quotient determines the Hopf subalgebra as the space of coinvariants (cf.~\cite{M}, Lemma~3.4.2 and Proposition~3.4.3). Therefore, $N'$~is uniquely determined, in other words, we have~$N'=N$. We have already seen in Paragraph~\ref{N16} that~$N$ is neither commutative nor cocommutative, so the corresponding extension is not abelian.
\qed
\end{pf}

As we mentioned in the introduction, many of the known semisimple Hopf algebras of prime power dimension contain a large commutative Hopf subalgebra of prime index. The result above shows that~$B$ is not of this kind.

\subsection[Hopf subalgebras]{} \label{HSA}
We have seen in Proposition~\ref{QN16} that~$N$ has exactly three quotient Hopf algebras of dimension~$8$. It has also exactly three Hopf subalgebras of dimension~$8$:
\begin{prop}
The Hopf algebra~$N$ has exactly three Hopf subalgebras of dimension~$8$, namely
\begin{align*}
M_1 \deq K\langle c_4, h_2, h_3 \rangle \qquad \quad
M_2 \deq K\langle c_2, c_4, h_3 \rangle \qquad \quad
M_3 \deq K\langle c_2 h_2, h_3 \rangle
\end{align*}
All these three subalgebras are commutative and normal in~$N$.
\end{prop}
\begin{pf}
\begin{pflist}
\item
Suppose that~$M$ is an eight-dimensional Hopf subalgebra of~$N$. As already mentioned twice, we then know from~\cite{KoMa}, Proposition~2 that~$M$ is normal in~$N$. The quotient~$N/NM^+$ has dimension~$2$, and the transpose of the quotient map yields an embedding
\[(N/NM^+)^* \to N^*\]
As in the proof of Theorem~\ref{Proof2by16}, the two-dimensional image of this map must have the form
$\Span(\varepsilon, \chi)$ for a nontrivial central group-like element~$\chi$ of order~$2$. This means that the original quotient map is essentially the map
\[N \to K \times K,~b \mapsto (\varepsilon(b), \chi(b)) \]
By~\cite{M}, Proposition~3.4.3, the original Hopf subalgebra can be recovered as the space of coinvariants, i.e., as the space
\[M = \{b \in N \mid b = b_\1 \chi(b_\2) \}\]
By tracing this argument backwards, we also see that, for a given nontrivial central group-like element~$\chi$ of order~$2$, the last equation yields an eight-dimensional normal Hopf subalgebra~$M$ of~$N$.

\item
From Proposition~\ref{N16}, we know that $N^*$~contains exactly three nontrivial central group-like elements, namely $\chi'_1$, $\chi'_2$, and~$\chi'_1 \chi'_2$. The computations in the proof of that proposition show that, if $b = a \star h \in N$, we have
\[b_\1 \chi'_1 (b_\2) = a_\1 \varepsilon'_2(a_\2) \star  h\]
Now the map $C \to C,~a \mapsto a_\1 \varepsilon'_2(a_\2)$
maps~$\omega_1$ and~$\omega_4$ to themselves and~$\omega_2$ and~$\omega_3$ to their negatives. This implies that
\[\{b \in N \mid b = b_\1 \chi'_1(b_\2) \}= K\langle c_4, h_2, h_3 \rangle = M_1\]

\item
For $b = a \star h \in N$, the computations in the proof of Proposition~\ref{N16} also show that
\[b_\1 \chi'_2 (b_\2) = a \star  h_\1 \gamma_3(h_\2) \]
Since $\gamma_3(g_3) = \theta(g_3, g_3) = 1$, but $\gamma_3(g_2) = \theta(g_3, g_2) = -1$, we see that
\[\{b \in N \mid b = b_\1 \chi'_2(b_\2) \} = K\langle c_2, c_4, h_3 \rangle = M_2\]

\item
Finally, in the case of the third element~$\chi'_1 \chi'_2$, the computations in the proof of Proposition~\ref{N16} yield for $b = a \star h \in N$ that
\[b_\1 (\chi'_1 \chi'_2)(b_\2) = a_\1 \varepsilon'_2(a_\2) \star  h_\1 \gamma_3(h_\2) \]
By combining the considerations for the first two elements, we see that
\begin{align*}
\{b \in N \mid b = b_\1 (\chi'_1 \chi'_2)(b_\2) \} &=
\Span(1, c_4, h_3, c_4 h_3, c_2 h_2, c_3 h_2, c_2 h_4, c_3 h_4)
\end{align*}
Now we know from Equation~(\ref{c2h2}) that $(c_2 h_2)^2 = c_4$, so that
$(c_2 h_2)^3 = c_3 h_2$ and $(c_2 h_2)^4 = 1$. This shows that this space is indeed equal to 
$M_3 = K\langle c_2 h_2, h_3 \rangle$.
\qed
\end{pflist}
\end{pf}

The first two of these eight-dimensional Hopf subalgebras have a very natural interpretation: According to Lemma~\ref{N16}, the subalgebra $M_1$~is just the span of the group-like
elements~${\bf G}(N)$. The Hopf algebra~$M_2$ is the biproduct
\[M_2 = C \star K[\langle g_3 \rangle] \]
which can be formed because the formulas for~$\delta(x)$ and~$\delta(x^3)$ in Paragraph~\ref{YetDrinfFirst} show that~$C$ is indeed a Yetter-Drinfel'd Hopf algebra
over~$K[\langle g_3 \rangle]$. Note that, because the action of~$g_3$ on~$C$ is trivial, $M_2$~is in fact the tensor product of~$C$ and~$K[\langle g_3 \rangle]$ as an algebra. The Hopf subalgebra~$M_3$ is less obvious. However, we will see in Paragraph~\ref{DHSA} that~$M_2$ and~$M_3$ are isomorphic.

The first two Hopf subalgebras can also be easily expressed in terms of the generators~$u$, $r$, and~$s$ introduced in Paragraph~\ref{AlgBiprod}, namely as
\begin{align*}
M_1 = K\langle u^2, r, s \rangle \qquad \qquad M_2 = K\langle u, s \rangle
\end{align*}
The third one looks more complicated when expressed in this way: $M_3$ is generated by $\frac{\K}{2}(1 + \iota \zeta^2) u r + \frac{\K}{2}(1 - \iota \zeta^2) u^3 r$ together with~$s$.

\subsection[Determination of the Hopf subalgebras]{} \label{DHSA}
Semisimple Hopf algebras of dimension~$8$ have been classified in~\cite{Ma1}. So the question arises which of these Hopf algebras the algebras~$M_1$, $M_2$, and~$M_3$ actually are. For~$M_1$, this is not difficult: As we already said above, it is the group algebra of the group~${\bf G}(N) \cong \Z_2 \times \Z_2 \times \Z_2$. The remaining two, being commutative, must be dual group algebras, as we mentioned in Paragraph~\ref{Extens}. It turns out that they are both isomorphic to the dual group algebra of the dihedral group~$D_8$:
\begin{prop}
We have $M_2 \cong M_3 \cong K^{D_8}$.
\end{prop}
\begin{pf}
\begin{pflist}
\item
Because $M_2 = K\langle c_2, c_4, h_3 \rangle$ is commutative and~$K$ is algebraically closed, we must have $\left|{\bf G}(M_2^*)\right|=8$. Since $c_2^2 = c_4^2 = h_3^2 = 1$, any multiplicative character~$\chi \in {\bf G}(M_2^*)$ satisfies
\[\chi (c_2) = \pm 1 \qquad \qquad \chi (c_4) = \pm 1 \qquad \qquad \chi (h_3) = \pm 1\]
and therefore, since  $\left|{\bf G}(M_2^*)\right|=8$, all possible combinations of $\pm 1$ must yield multiplicative characters.

For two multiplicative characters $\chi, \chi' \in {\bf G}(M_2^*)$, we have on the last two generators 
\[(\chi \chi')(c_4) = \chi(c_4) \chi'(c_4) \qquad \text{and} \qquad
(\chi \chi')(h_3) = \chi(h_3) \chi'(h_3)\]
because these generators are group-like. For the third generator~$c_2$, the formula
\begin{align*}
\db(c_2) &= \frac{1}{2} c_2 (\B + h_3) \ot c_2 + \frac{1}{2} c_2 (\B - h_3) \ot c_2c_4
\end{align*}
for its coproduct, which was given in Paragraph~\ref{CoalgBiprod} and in particular shows that~$M_2$ is not cocommutative, implies that
\begin{align*}
(\chi \chi') \left( c_2 \right)  &=
\frac{1}{2} \left(1 + \chi(h_3) + \chi'(c_4) - \chi(h_3) \chi'(c_4)\right) \chi(c_2) \chi'(c_2) \\
&= \begin{cases}
\chi(c_2) \chi'(c_2) & : \chi(h_3) = 1
\mspace{10mu} \text{or} \mspace{10mu} \chi'(c_4) = 1\\
- \chi(c_2) \chi'(c_2) & : \chi(h_3) = \chi'(c_4) = -1
\end{cases}
\end{align*}
If we apply this formula in the case $\chi = \chi'$, it shows that ${\bf G}(M_2^*)$ has exactly two elements of order~$4$, defined by $\chi(c_2) = \pm 1$ and $\chi(c_4) = \chi(h_3) =- 1$, and five elements of order~$2$, which can only happen if~\mbox{${\bf G}(M_2^*) \cong D_8$}.

\item
In $M_3 = K \langle c_2 h_2, h_3 \rangle$, we have by Equation~(\ref{c2h2}) that
$(c_2 h_2)^2 = c_4$, so that $(c_2 h_2)^3 = c_3 h_2$ and $(c_2 h_2)^4 = 1$, as we have already recalled in the proof of Proposition~\ref{HSA}. Each of the eight multiplicative characters therefore maps the first generator~$c_2 h_2$ to a fourth root of unity and the second generator~$h_3$ to~$\pm 1$, and any fourth root of unity and any sign can arise in this way. The formula for the coproduct of~$c_2$ given in Paragraph~\ref{CoalgBiprod} now yields
\begin{align*}
\db(c_2 h_2) &=
\frac{1}{2}
\left(( c_2 h_2 + c_2 h_2 h_3) \ot c_2 h_2 + (c_2 h_2 - c_2 h_2 h_3) \ot c_3 h_2\right) \\
&= \frac{1}{2}
\left(( c_2 h_2 + c_2 h_2 h_3) \ot c_2 h_2 + (c_2 h_2 - c_2 h_2 h_3) \ot (c_2 h_2)^3\right)
\end{align*}
This shows on the one hand that~$M_3$ is a subbialgebra, and therefore a Hopf subalgebra, of~$B$, and on the other hand that~$M_3$ is not cocommutative.

For two multiplicative characters $\chi, \chi' \in {\bf G}(M_3^*)$, we have for the second generator that
$(\chi \chi')(h_3) = \chi(h_3) \chi'(h_3)$
and for the first generator that
\begin{align*}
(\chi \chi')(c_2 h_2)
&= \frac{1}{2} \left( 1 + \chi(h_3)  + \chi'(c_2 h_2)^2 - \chi(h_3) \chi'(c_2 h_2)^2 \right)
\chi(c_2 h_2) \chi'(c_2 h_2)\\
&= \begin{cases}
\chi(c_2 h_2) \chi'(c_2 h_2) & : \chi(h_3) = 1 \\
\chi(c_2 h_2) \chi'(c_2 h_2)^3 & : \chi(h_3) = -1
\end{cases}
\end{align*}
If we apply this formula in the case $\chi = \chi'$, it shows that~${\bf G}(M_3^*)$ has also exactly two elements of order~$4$, defined by $\chi(c_2 h_2) = \pm \iota$ and $\chi(h_3) = 1$, and five elements of order~$2$, so that again~\mbox{${\bf G}(M_3^*) \cong D_8$}.
\qed
\end{pflist}
\end{pf}

We note that there is a second way to see that~$N$ has exactly three eight-dimensional Hopf subalgebras: It follows from the classification of semisimple Hopf algebras of dimension~$16$ that $N \cong H_{d:1,1}$ is self-dual (cf.~\cite{Ka}, Remark~1.4). Since~$N$ has, by Proposition~\ref{QN16}, up to isomorphism exactly three quotient Hopf algebras of dimension~$8$, self-duality implies that~$N \cong N^*$ has exactly three Hopf subalgebras of dimension~$8$. Our discussion above, however, gives their explicit form.

\section{Hopf subalgebras of dimension~$2$} \label{Sec:Dim2}
\subsection[Quotients]{} \label{YDQ}
The four-dimensional Yetter-Drinfel'd Hopf subalgebra~$C$ of~$A$ introduced in Paragraph~\ref{YetDrinfFirst} also arises as a quotient of~$A$:
\begin{prop}
There is a unique homomorphism $\pi_C \colon A \to C$ of Yetter-Drinfel'd Hopf algebras defined on generators as
\[\pi_C(x) \deq \omega_4 = x^2 \qquad \qquad \pi_C(y) \deq \omega_2 \]
Its values on the group-like elements are given by the following table:
\begin{center}
\renewcommand\arraystretch{1.5}
\begin{tabular}{|c|c|c|c|c|c|c|c|c|} \hline
$a$ & $\omega_1$ & $\omega_2$ & $\omega_3$ & $\omega_4$
& $\eta_1$ & $\eta_2$ & $\eta_3$ & $\eta_4$\\ \hline
$\pi_C(a)$ & $\omega_1$ & $\omega_4$ & $\omega_4$ & $\omega_1$
& $\omega_2$ & $\omega_3$ & $\omega_2$ & $\omega_3$\\ \hline
\end{tabular}
\end{center}
\end{prop}
\begin{pf}
It is immediate from the defining relations of~$A$ that there is such an algebra homomorphism, and it is also immediate from the definitions that it takes the specified values on the group-like elements. Because it takes group-like elements to group-like elements, it is a coalgebra homomorphism that satisfies~$\sd \circ \pi_C = \pi_C \circ \sa$. It commutes with the action of~$g_2$ since
\[\pi_C(g_2.x) = \pi_C(x^3) = x^2 = g_2.\pi_C(x)\]
and
\[\pi_C(g_2.y) = \pi_C(x^3 y) = \omega_4 \omega_2 = \omega_3 = g_2.\pi_C(y)\]
Similarly, it commutes with the action of~$g_3$ since
\[\pi_C(g_3.x) = \pi_C(x) = x^2 = g_3.\pi_C(x)\]
and
\[\pi_C(g_3.y) = \pi_C(x^2 y) = \omega_2 = g_3.\pi_C(y)\]
It is therefore $H$-linear and, because the coaction is derived from the action via Equation~(\ref{coaction}), also colinear.
\qed
\end{pf}

On the Radford biproduct, $\pi_C$ induces the surjective Hopf algebra homomorphism
\[\pi_N \colon B \to N,~a \star h \mapsto \pi_C(a) \star h \]
which on the generators introduced in Paragraph~\ref{AlgBiprod} takes the values
\[\pi_N(u) = u^2 \qquad \qquad
\pi_N(v) = \frac{\K}{2}(1 + \iota \zeta^2) u + \frac{\K}{2}(1 - \iota \zeta^2) u^3  \]
as well as $\pi_N(r) = r$ and~$\pi_N(s) = s$. Note that the second equation can also be stated as~$\pi_N(d_1) = c_2$. By comparing the table in Paragraph~\ref{Group} to the table in Paragraph~\ref{N16}, we see that~$\pi_N$ relates the group-like elements of~$B^*$ and the group-like elements of~$N^*$ via
\[\chi_1 = \chi'_1 \circ \pi_N \qquad  \qquad \chi_2 = \chi'_2 \circ \pi_N \qquad \qquad
\chi_3 = \chi'_3 \circ \pi_N \qquad \qquad\]
We will denote the kernel of~$\pi_N$ by~$I$, so that~$B/I \cong N$.

\subsection[Uniqueness of quotients]{} \label{16by2}
The ideal~$I$ just defined is a Hopf ideal of dimension~$16$. We will now show that it is the only such ideal. In preparation, we prove a simple lemma:
\begin{lem}
The space~$U \deq \Span(c_1, c_4)$ is the only normal Hopf subalgebra of~$B$ of dimension~$2$. It is central.
\end{lem}
\begin{pf}
We have already stated in Paragraph~\ref{Extens} that a normal Hopf subalgebra~$U'$ of~$B$ of dimension~$2$ is central. Since it is commutative, cocommutative, and semisimple, $U'$~is spanned by the unit element and a nontrivial central group-like element. But according to Lemma~\ref{Group}, the only nontrivial central group-like element is~$c_4$, and so the assertion holds.
\qed
\end{pf}

The Hopf ideal~$I$ and the Hopf subalgebra~$U$ are related as follows:
\begin{thm}
The ideal~$I$ is the only Hopf ideal of~$B$ of dimension~$16$. It is normal. We have~$I = BU^+$ and $U = B^{\co N}$. The Hopf algebra~$B$ therefore fits into exactly one extension of the type
\[U \hookrightarrow B \overset{\pi_N}{\twoheadrightarrow} N\]
with~$\dm U = 2$ and~$\dm N = 16$. Since~$N$ is neither commutative nor cocommutative, this extension is not abelian.
\end{thm}
\begin{pf}
If~$I'$ is a Hopf ideal of~$B$ of dimension~$16$, we can set~$N' \deq B/I'$ and denote the quotient map by~$\pi_{N'} \colon B \to N'$. Then $\pi_{N'}^*(N'^*)$ is a Hopf subalgebra of~$B^*$ of index~$2$. As in the proof of Theorem~\ref{Proof2by16}, it follows from~\cite{KoMa}, Proposition~2 that $\pi_{N'}^*(N'^*)$ is normal in~$B^*$. This implies that~$I'$ is normal and that~$\pi_{N'}$ is a conormal surjection (cf.~\cite{M}, Section~3.4, Page~36), which in turn implies that~$U' \deq B^{\co N'}$ is a normal Hopf subalgebra of dimension~$2$. By the preceding lemma, we have $U'=U$, and consequently~$I'=BU^+$. This shows that~$I'$ is uniquely determined, and therefore~$I'=I$.
\qed
\end{pf}

\subsection[The Grothendieck ring of~$B$]{} \label{GrothB}
As promised in Paragraph~\ref{Wedderb}, the preceding results enable us to describe the product in the Grothendieck ring~$K_{0}(B)$. From Proposition~\ref{GrouplB*}, we know that
\[{\bf G}(B^*) = \left\langle \chi_1 \right\rangle \times \left\langle \chi_2 \right\rangle \times \left\langle \chi_3 \right\rangle \cong \Z_2 \times \Z_2 \times \Z_2\]
which describes the one-dimensional representations. By Theorem~\ref{RepsBB*}, the Hopf algebra~$B$ has two irreducible representations of dimension~$2$, which we will denote by~$\pi_{1}$ and~$\pi_{2}$, and one irreducible representation of dimension~$4$, which we will denote by~$\rho$.

In Theorem~\ref{16by2}, we have proved that~$B$ has a unique quotient Hopf algebra~$N$ of dimension~$16$, which is isomorphic to the Hopf algebra~$H_{d:1,1}$ from~\cite{Ka}, Table~1. Clearly, the irreducible representations of~$N$ can be pulled back along the quotient map~$\pi_N$ to irreducible representations of~$B$, so that we obtain a ring homomorphism from~$K_{0}(N)$ to~$K_{0}(B)$. The one-dimensional and the two-dimensional irreducible representations of~$B$ arise in this way; for the one-dimensional representations, we have already seen that at the end of Paragraph~\ref{YDQ}, and we have also indicated at the end of Paragraph~\ref{H16} how the one-dimensional representations of~$N$ correspond to the one-dimensional representations of~$H_{d:1,1}$.

Therefore, the description of the Grothendieck ring of~$H_{d:1,1}$ given in~\cite{Ka}, Paragraph~5.1 implies on the one hand that the two-dimensional irreducible representations, and also the one-dimensional representations, are self-dual, and on the other hand that the products of the one-dimensional with the two-dimensional irreducible representations are given by the table
\begin{center}
\renewcommand\arraystretch{1.5}
\begin{tabular}{|c|c|c|c|} \hline
& $\chi_1$ & $\chi_2$ & $\chi_3$  \\ \hline
$\pi_1$ & $\pi_1$ & $\pi_1$ & $\pi_2$ \\ \hline
$\pi_2$ & $\pi_2$ & $\pi_2$ & $\pi_1$ \\ \hline
\end{tabular}
\end{center}
while the products of the two-dimensional irreducible representations are given by the table
\begin{center}
\renewcommand\arraystretch{1.5}
\begin{tabular}{|c|c|c|} \hline
& $\pi_1$ & $\pi_2$ \\ \hline
$\pi_1$ & $1 + \chi_1 + \chi_2 + \chi_1 \chi_2$
& $\chi_3 + \chi_1 \chi_3 + \chi_2 \chi_3 + \chi_1 \chi_2 \chi_3$ \\ \hline
$\pi_2$ & $\chi_3 + \chi_1 \chi_3 + \chi_2 \chi_3 + \chi_1 \chi_2 \chi_3$
& $1 + \chi_1 + \chi_2 + \chi_1 \chi_2$ \\ \hline
\end{tabular}
\end{center}
Note that the self-duality of the irreducible representations implies that the Grothendieck ring of~$N$ is commutative, so that, in the above tables, it does not matter which factor in the product comes first.

It remains to treat the products that involve~$\rho$. Because it is the only four-dimensional irreducible representation, it is also self-dual, and we have \mbox{$\chi_i \rho = \rho$}
for~$i=1, 2, 3$. Furthermore, as in the case of~$N$, self-duality implies that the Grothendieck ring of~$B$ is commutative.

Since~$\rho$ does not occur in the decomposition of~$\chi_i \pi_j$ and~$\pi_i \pi_j$, it follows from~\cite{NR}, Theorem~9 that no one-dimensional or two-dimensional irreducible representation  occurs in the decomposition of~$\pi_j \rho$, so that
\[\pi_1 \rho = \pi_2 \rho = 2 \rho\]
Similarly, since~$\rho$ occurs in the decomposition of~$\chi_i \rho$ with multiplicity~$1$ and in the decomposition of~$\pi_i \rho$ with multiplicity~$2$, we get in the same way from \cite{NR}, Theorem~9 that
\[\rho^2 = 1 + \chi_1 + \chi_2 +\chi_3 + \chi_1 \chi_2 + \chi_1 \chi_3 + \chi_2 \chi_3
+ \chi_1 \chi_2 \chi_3 + 2 \pi_1 +2 \pi_2 \]
Note that the first eight summands are precisely the eight elements of~${\bf G}(B^*)$.

\section{Hopf subalgebras of dimensions~$4$ and~$8$} \label{Sec:4by8}
\subsection[The relation between~$B$ and~$N$]{} \label{BN}
From Theorem~\ref{Proof2by16} and Theorem~\ref{16by2}, we know that neither~$B$ nor~$B^*$ fits into an abelian extension of a Hopf algebra of dimension~$2$ by a commutative Hopf algebra of dimension~$2^n$. Thus neither~$B$ nor~$B^*$ can be constructed as an extension of the type used in~\cite{Ka}, \cite{Ka2}, \cite{Ka3}, and~\cite{Ma1} to classify certain semisimple Hopf algebras of prime power dimension. We will now show that nevertheless~$B$ fits into an abelian extension of  a Hopf algebra of dimension~$4$ by a commutative Hopf algebra of dimension~$8$, but this extension, which is essentially unique, is neither central nor cocentral.
\begin{thm} \label{all8}
The Hopf algebra~$B$ has exactly three Hopf subalgebras of dimension~$8$, namely
\begin{align*}
M_1 \deq K\langle c_4, h_2, h_3 \rangle \qquad \quad
M_2 \deq K\langle c_2, c_4, h_3 \rangle \qquad \quad
M_3 \deq K\langle c_2 h_2, h_3 \rangle
\end{align*}
\end{thm}
\begin{pf}
Suppose that~$M$ is an eight-dimensional Hopf subalgebra of~$B$. From Theorem~\ref{RepsBB*}, we know that~$B^*$ has eight $1$-dimensional, two $2$-dimensional, and one $4$-dimensional irreducible representations. Furthermore, it follows from Theorem~\ref{Proof2by16} that the unique quotient Hopf algebra of~$B^*$ of dimension~$16$ is the dual~$N^*$ of the Hopf subalgebra~$N$ of~$B$ defined in Paragraph~\ref{N16}. Since no Hopf algebra of dimension~$16$ or~$8$ can have a four-dimensional irreducible representation, the canonical restriction mappings from~$B^*$ to~$N^*$ and~$M^*$, which are dual to the inclusion mappings, must both contain the sixteen-dimensional ideal corresponding to the four-dimensional irreducible representation in their respective kernels. For the canonical map from~$B^*$ to~$N^*$, this kernel must in fact be equal to this ideal. This shows that the restriction mapping from~$B^*$ to~$M^*$ factors over the restriction mapping from~$B^*$ to~$N^*$, which implies that~$M \subset N$. But then Proposition~\ref{HSA} implies that~$M=M_1$, $M=M_2$, or~$M=M_3$.
\qed
\end{pf}

We have seen in Proposition~\ref{HSA} that~$M_1$, $M_2$, and~$M_3$ are normal in~$N$. They are not all normal in~$B$:
\begin{prop} \label{notnormal8}
$M_1$ and~$M_3$ are not normal in~$B$.
\end{prop}
\begin{pf}
\begin{pflist}
\item
If~$M_1$ were normal in~$B$, the ideal $B M_1^+ = M_1^+ B$ would be two-sided, and in the corresponding quotient algebra $B/(B M_1^+)$, which would be four-dimensional by~\cite{M}, Corollary~8.4.7, we would have
\[\bar{c}_4 = \bar{h}_2 = \bar{h}_3 = 1\]
Therefore we would also have $\bar{h}_4 = \bar{h}_2 \bar{h}_3 = 1$ and
$\bar{c}_3 = \bar{c}_4 \bar{c}_2 = \bar{c}_2$ as well as
$\bar{d}_3 = \bar{c}_4 \bar{d}_1 = \bar{d}_1$ and $\bar{d}_4 = \bar{c}_4 \bar{d}_2 = \bar{d}_2$, so that
\[B/(B M_1^+) = \Span(1, \bar{c}_2, \bar{d}_1, \bar{d}_2)\]
But as we would then also have
\[\bar{d}_1 = \bar{h}_2 \bar{d}_1 = \bar{d}_2 \bar{h}_2 = \bar{d}_2\]
it would follow that $B/(B M_1^+) = \Span(1, \bar{c}_2, \bar{d}_1)$, which would be a contradiction.

\item
The argument showing that~$M_3$ is not a normal Hopf subalgebra of~$B$ is very similar. If it were, $BM_3^+ = M_3^+ B$ would be a two-sided ideal of~$B$, and the corresponding quotient algebra $B/(B M_3^+)$ would again have dimension~$4$. By considering the basis of~$M_3$ given in the proof of Proposition~\ref{HSA}, we see that we would have
\[\bar{c}_4 = \bar{h}_3 = \bar{c}_4 \bar{h}_3 = \bar{c}_2 \bar{h}_2 = \bar{c}_3 \bar{h}_2
= \bar{c}_2 \bar{h}_4 = \bar{c}_3 \bar{h}_4 = 1\]
in the quotient algebra. In this quotient, we would also have
\[\bar{c}_2 = \bar{c}_4 \bar{c}_2 = \bar{c}_3= \bar{h}_4 = \bar{h}_3 \bar{h}_2 = \bar{h}_2\]
as well as
$\bar{d}_3 = \bar{c}_4 \bar{d}_1 = \bar{d}_1$ and $\bar{d}_4 = \bar{c}_4 \bar{d}_2 = \bar{d}_2$, so that
\[B/(B M_3^+) = \Span(1, \bar{c}_2, \bar{d}_1, \bar{d}_2)\]
In this case, we would have
\[\bar{h}_2 \bar{d}_1 = \bar{d}_2 \bar{h}_2 = \bar{d}_2 \bar{c}_2 = \bar{c}_2 \bar{d}_2 
= \bar{h}_2 \bar{d}_2\]
and therefore $\bar{d}_1 = \bar{d}_2$, so that we would arrive again at the contradiction 
that~$B/(B M_3^+) = \Span(1, \bar{c}_2, \bar{d}_1)$.
\qed
\end{pflist}
\end{pf}

\subsection[Normal Hopf subalgebras of dimension~$8$]{} \label{normal8}
On the other hand, there is a normal eight-dimensional Hopf subalgebra:
\begin{lem}
$M_2$ is normal in~$B$.
\end{lem}
\begin{pf}
Recall that~$M_2$ is generated by~$c_2$, $c_4$, and~$h_3$. From the description of the generators $\chi_1, \chi_2, \chi_3$ of~${\bf G}(B^*)$ in Proposition~\ref{Group}, we know that
\[\chi_1(c_2) = \chi_1(c_4) = \chi_1(h_3) = \K\]
as well as $\chi_2(c_2) = \chi_2(c_4) = \chi_2(h_3) = \K$, so that $\chi_1$ and~$\chi_2$ take the generators of~$M_2$ to~$\K$. By Corollary~\ref{N16} and the remark following it, this also holds for~$\chi_1 \chi_2$.

As in the proof of Proposition~\ref{HSA}, the algebra homomorphism
\[\pi \colon B \to K^4,~b \mapsto (\eb(b), \chi_1(b), \chi_2(b), (\chi_1 \chi_2)(b))\]
is essentially the transpose of the inclusion map from
$\Span(\eb, \chi_1, \chi_2, \chi_1 \chi_2)$ to~$B^*$, and is therefore surjective. Because~$c_4$ and~$h_3$ are group-like, our computations above show that these two generators are contained in the spaces~$B^{\co \pi}$ and~$\prescript{\co \pi}{} \!B$ of right and left coinvariants, which are eight-dimensional by~\cite{M}, Corollary~8.4.7. The formula for the coproduct of~$c_2$ given in Paragraph~\ref{CoalgBiprod} shows that this also holds for~$c_2$, so that~$M_2$ is completely contained in these spaces. But as~$M_2$ is itself eight-dimensional, we have 
\[B^{\co \pi} = \prescript{\co \pi}{} \!B = M_2\]
By~\cite{M}, Lemma~3.4.2, this implies that~$M_2$ is normal in~$B$.
\qed
\end{pf}

The preceding proof also yields the form of the quotient:
\begin{thm}
The Hopf subalgebra~$M \deq M_2$ is the unique normal Hopf subalgebra of~$B$ of dimension~$8$. Thus~$B$ fits into exactly one extension of the type
\[M \hookrightarrow B \twoheadrightarrow Q\]
with $\dm M = 8$ and $\dm Q = 4$. In this case, $M \cong K^{D_8}$
and~$Q \cong K[\Z_2 \times \Z_2]$. This extension is abelian, but neither central nor cocentral.
\end{thm}
\begin{pf}
\begin{pflist}
\item
That~$M_2$ is the unique normal Hopf subalgebra of~$B$ of dimension~$8$ follows from Theorem~\ref{BN}, Proposition~\ref{BN}, and the lemma above. Therefore, we have a unique extension with the stated properties. From Proposition~\ref{DHSA}, we know
that~$M_2 \cong K^{D_8}$.

\item
The proof of the lemma above also yields that the quotient \mbox{$Q \deq B/(B M_2^+)$} is isomorphic to the dual of
$\Span(\eb, \chi_1, \chi_2, \chi_1 \chi_2) \cong K[\Z_2 \times \Z_2]$. Since this Hopf algebra is self-dual, we have \mbox{$Q \cong K[\Z_2 \times \Z_2]$} as well. Without referring to the proof of the lemma, this can be seen as follows: Since $\dm Q = 4$, the quotient must be a group algebra~$Q=K[L]$ for a group~$L$ of order~$4$, so that~$L \cong \Z_2 \times \Z_2$ or~$L \cong \Z_4$. Then~$Q^*$ is isomorphic to the group algebra~$K[\hat{L}]$ of the character group, and the transpose~$\pi_Q^*$ of the quotient map~$\pi_Q$ yields an injective group homomorphism from~$\hat{L}$ to~${\bf G}(B^*)$, which is elementary abelian of order~$8$ by Proposition~\ref{Group}. Therefore $L \cong \hat{L}$ must be elementary abelian of order~$4$.

\item
The extension is not central, for example because~$M_2$ contains the noncentral element~$c_2$.
From Proposition~\ref{Group}, we know that~$B^*$ has only one nontrivial central group-like element. So~$\pi_Q^*(Q^*)$ cannot be central in~$B^*$, which means that this extension is not cocentral.
\qed
\end{pflist}
\end{pf}

\subsection[Normal Hopf subalgebras of dimension~$4$]{} \label{normal4}
We have just seen that~$B$ fits into exactly one abelian extension of a Hopf algebra of dimension~$4$ by a Hopf algebra of dimension~$8$. We will now show that~$B$ also fits into exactly one abelian extension of a Hopf algebra of dimension~$8$ by a Hopf algebra of dimension~$4$, which is also neither central nor cocentral.
\begin{thm} \label{all4}
The Hopf algebra~$B$ has seven Hopf subalgebras of dimension~$4$. Precisely one of them, namely $K\langle c_4, h_3 \rangle$, is normal.
\end{thm}
\begin{pf}
\begin{pflist}
\item Since every semisimple Hopf algebra of dimension~$4$ over an algebraically closed field of characteristic zero is isomorphic to a group algebra, every Hopf subalgebra of~$B$ of dimension~$4$ is a group algebra of a subgroup of order~$4$ of~${\bf G}(B)$. From Lemma~\ref{Group}, we know that~${\bf G}(B)$ is an elementary abelian group of order~$8$, or equivalently a three-dimensional vector space over the field with two elements. A subgroup of order~$4$ is then the same as a two-dimensional subspace. This subspace has a one-dimensional orthogonal complement in the dual vector space, and this orthogonal complement determines the subspace. A~one-dimensional subspace over the field with two elements is determined by its unique nonzero vector. In the dual space, there are seven nonzero vectors, and correspondingly there are seven subgroups of order~$4$ of~${\bf G}(B)$. Explicitly, these are the groups $G_1 \deq \langle c_4, h_3 \rangle = \{\B, c_4, h_3, c_4 h_3\}$ and
\begin{align*}
G_2 &\deq \{\B, c_4, h_2, c_4 h_2\} & 
G_3 &\deq \{\B, c_4, h_4, c_4 h_4\} \\
G_4 &\deq \{\B, h_2, h_3, h_4\} & G_5 &= \{\B, h_2, c_4 h_3, c_4 h_4\} \\
G_6 &\deq \{\B, h_3, c_4 h_2, c_4 h_4\} &
G_7 &\deq \{\B, h_4, c_4 h_2, c_4 h_3\}
\end{align*}

\item
Suppose that, for one of these groups~$G_i$, the group algebra~$P \deq K[G_i]$ is normal in~$B$. Because~${\bf G}(B) = {\bf G}(N)$, we then have~$P \subset N$. It follows directly from~\cite{M}, Definition~3.4.1 that~$P$ is also normal in~$N$. By~\cite{M}, Corollary~8.4.7, the quotient~$N/(NP^+)$ has dimension~$4$. It is therefore a group algebra of an abelian group of order~$4$. The four distinct multiplicative characters
in~${\bf G}((N/(NP^+))^*)$ yield by pullback along the quotient map a subgroup of order~$4$ in~${\bf G}(N^*)$ whose elements vanish on~$NP^+$, and therefore take the value~$\K$ on~$G_i$.

Now we see from Corollary~\ref{N16} that the subgroup of~${\bf G}(N^*)$ of characters that take the value~$\K$ on the group generated by~$r=h_2$ and~$s=h_3$ has only order~$2$ and is generated by~$\chi'_1$. On the other hand, every element of~${\bf G}(N^*)$ takes the value~$\K$ on~$u^2=c_4$. Thus the fact that there are four distinct elements of~${\bf G}(N^*)$ that map~$G_i$ to~$\K$ rules out that $i=4$, $5$, $6$, or~$7$; in other words,~$G_i$ must contain~$c_4$, and therefore~$P$ must contain~$U \deq \Span(c_1, c_4)$.

\item
Consequently, we have a surjective map~$B/(BU^+) \to B/(BP^+)$. By composing it with the inverse of the isomorphism $\bar{\pi}_N \colon B/(BU^+) \to N$ from Theorem~\ref{16by2}, we obtain an eight-dimensional quotient of~$N$. From Proposition~\ref{QN16}, we know that~$N$~has exactly three quotients of dimension~$8$, denoted there by~$F_1$, $F_2$, and~$F_3$.

The multiplicative character~$\chi'_2 \in {\bf G}(N^*)$ satisfies
$\chi'_2(c_4) = \chi'_2(u^2) = \K$ and $\chi'_2(h_3) = \chi'_2(s) = \K$ and therefore factors to a multiplicative character on all three quotients~$F_1$, $F_2$, and~$F_3$. As discussed in Paragraph~\ref{YDQ}, we have \mbox{$\chi_2 = \chi'_2 \circ \pi_N$}, so~$\chi_2$ factors to a multiplicative character of~$B/(BP^+)$. In particular, $\chi_2$ maps all elements of~$G_i$ to~$\K$. This implies that~$i=1$. We have therefore shown that, if~$P$ is a four-dimensional normal Hopf subalgebra of~$B$, then~$P=K[G_1]$.

\item
\enlargethispage{4pt}
To see that~$K[G_1]$ is indeed normal in~$B$, we argue as in the proof of Lemma~\ref{normal8}: Recall that $F_2 = N/\left( K \langle h_3 \rangle^+ N \right)$ and consider the Hopf algebra homomorphism that arises as the composition
\[B \overset{\pi_N}{\to} N \to F_2\]
This composition maps the group-like elements~$c_4$ and~$h_3$ to~$1_{F_2}$, so that they are both contained in the spaces~$B^{\co F_2}$ and~$\prescript{\co F_2}{} \!B$ of right and left coinvariants. As these spaces are four-dimensional by~\cite{M}, Corollary~8.4.7, we have 
\[B^{\co F_2} = \prescript{\co F_2}{} \!B = K[G_1]\]
By~\cite{M}, Lemma~3.4.2, this implies that~$K[G_1]$ is normal in~$B$.
\qed
\end{pflist}
\end{pf}

There is a somewhat more direct way to see that~$K[G_i]$ is not normal in~$B$ for~$i \neq 1$, which uses an argument that is similar to the proof of Proposition~\ref{BN} and proceeds on a case by case basis. We illustrate this argument by giving it in the case~$i=3$. If we assume that~$P=K[G_3]$ is normal in~$B$, we know from~\cite{M}, Corollary~3.4.4 that~$BP^+= P^+B$. In the eight-dimensional quotient algebra~$B/(BP^+)$, we have 
\[\bar{c}_4 = \bar{h}_4 = 1\]
and therefore also $\bar{h}_2 = \bar{h}_3$, $\bar{c}_2 = \bar{c}_3$, $\bar{d}_1 = \bar{d}_3$, and $\bar{d}_2 = \bar{d}_4$.
This implies that
\[B/(BP^+) = \Span(1, \bar{c}_2, \bar{d}_1, \bar{d}_2, \bar{h}_2, \bar{c}_2 \bar{h}_2, \bar{d}_1\bar{h}_2, \bar{d}_2 \bar{h}_2)\]
But since
\[\bar{d}_1 = \bar{h}_4 \bar{d}_1 = \bar{d}_4 \bar{h}_4 = \bar{d}_4 = \bar{d}_2\]
the algebra~$B/(BP^+)$ is spanned by only seven elements, which is a contradiction.

The preceding theorem has the following consequence:
\begin{cor}
The Hopf algebra~$B$ fits into exactly one extension of the type
\[P \hookrightarrow B \twoheadrightarrow F \]
with $\dm P = 4$ and $\dm F =8 $. In this case, $P = K \langle c_4, h_3 \rangle \cong K[\Z_2\times \Z_2]$ and $F \cong K[D_8]$. This extension is abelian, but neither central nor cocentral.
\end{cor}
\begin{pf}
We have seen in the preceding theorem that $P = K \langle c_4, h_3 \rangle$ is the unique normal Hopf subalgebra of~$B$ that has dimension~$4$, but in the proof we have also seen that the arising quotient is isomorphic to~$F_2$. From Paragraph~\ref{QN16}, we know that
$F_2 \cong K[D_8]$. The corresponding extension is therefore abelian. It is not a central extension, because~$h_3$ is not central.

The commutator factor group of~$D_8$ is isomorphic to~$\Z_2 \times \Z_2$, and therefore~$D_8$ has four one-dimensional representations. This implies
that~${\bf G}(F_2^*) \cong \Z_2 \times \Z_2$. But from Proposition~\ref{Group}, we know that~$B^*$ has only one nontrivial central group-like element, so $F_2^*$ cannot be contained in the center of~$B^*$. In other words, this extension is not cocentral.
\qed
\end{pf}

We note that, in terms of the generators introduced at the beginning of Paragraph~\ref{AlgBiprod}, we have $P = K \langle u^2, s \rangle$.

\section{The biproduct in the second case} \label{Sec:BiprodSec}
\subsection[The Yetter-Drinfel'd Hopf algebra in the second case]{} \label{YetDrinfSec}
So far, we have only treated the Yetter-Drinfel'd Hopf algebra~$A$ and its associated Radford biproduct~$B = A \star H$. In~\cite{KaSo2}, Section~3, the authors also described a second example of a Yetter-Drinfel'd Hopf algebra over~$H$ of dimension~$8$, denoted here by~$\Au$. As the algebra~$A$ introduced in Paragraph~\ref{YetDrinfFirst}, it is a Yetter-Drinfel'd Hopf algebra over the group algebra~$H=K[G]$ of the elementary abelian group~$G = \Z_2 \times \Z_2$, and it is also generated by two elements, denoted here by~$\xu$ and~$\yu$. But in contrast to the previous case, these generators do not commute and satisfy instead the defining relations
\[\xu^4 = 1 \qquad \qquad \xu \yu = \yu \xu^3 \qquad \qquad
\yu^2 = \frac{\K}{2}(\zeta 1 + \xu - \zeta \xu^2 + \xu^3)\]
where~$\zeta$ is again a not necessarily primitive fourth root of unity.

The action of~$G$ on~$\Au$ is defined by the same formulas as the action of~$G$ on~$A$, and the coaction is again derived from the action via Equation~(\ref{coaction}), using the same bicharacter~$\theta$. As in Paragraph~\ref{YetDrinfFirst}, we introduce the elements $\omu_1 \deq 1_{\Au}$,
\[\omu_2 \deq \frac{\K}{2}(1 + \iota \zeta^2) \xu + \frac{\K}{2}(1 - \iota \zeta^2) \xu^3 \qquad
\omu_3 \deq \frac{\K}{2}(1 - \iota \zeta^2) \xu + \frac{\K}{2}(1 + \iota \zeta^2) \xu^3\]
$\omu_4 \deq \xu^2$, and
\[\etu_1 \deq \yu \qquad \etu_2 \deq \xu^3 \yu \qquad \etu_3 \deq \xu^2 \yu \qquad \etu_4 \deq \xu \yu\]
where~$\iota$ is a fixed primitive fourth root of unity unrelated to~$\zeta$. According to \cite{KaSo2}, Proposition~3.4, these eight elements form a basis of~$A$, and the coalgebra structure is determined by the fact that they are group-like. By~\cite{KaSo2}, Paragraph~3.3, this means for~$\xu$ that
\[\Delta_{\Au}(\xu) =
\frac{\K}{2} (\xu \ot \xu + \xu \ot \xu^3 + \xu^3 \ot \xu - \xu^3 \ot \xu^3)\]
As in Paragraph~\ref{AlgBiprod}, we can form the biproduct~$\Bu \deq \Au \star H$. It is generated by the four elements
\[\uu \deq \xu \star \HH \qquad \quad \vu \deq \yu \star \HH \qquad \quad
\ru \deq 1_{\Au} \star g_2 \qquad \quad \su \deq 1_{\Au} \star g_3\]
In the same way as before, it is possible to derive a presentation of~$\Bu$:
\enlargethispage{4pt}
\begin{prop}
The four generators satisfy the relations
\begin{proplist}
\item
$\displaystyle
\uu^4 = 1, \qquad \uu \vu = \vu \uu^3, \qquad
\vu^2 = \frac{\K}{2}(\zeta 1 + \uu - \zeta \uu^2 + \uu^3)$

\item
$\displaystyle
\ru^2 = 1, \qquad \ru \su = \su \ru, \qquad \su^2 = 1$

\item
$\displaystyle
\ru \uu = \uu^3 \ru, \qquad \ru \vu = \uu^3 \vu \ru, \qquad \su \uu = \uu \su, \qquad
\su \vu = \uu^2 \vu \su$
\end{proplist}
These relations are defining.
\end{prop}

We can also, as in Paragraph~\ref{AlgBiprod}, introduce the elements that correspond to the basis elements of~$\Au$ by defining
\[\cu_i \deq \omu_i \star \HH \qquad \qquad \du_i \deq \etu_i \star \HH \qquad \qquad
\hu_i \deq 1_{\Au} \star g_i \]
for $i=1,2,3,4$. Then the elements $\cu_i \hu_j$ together with the elements $\du_i \hu_j$,
for~$i,j = 1,2,3,4$, form a basis of~$\Bu$. Note that $\cu_1 = \hu_1 = 1_{\Bu}$. Except for~$\uu$, the above generators are among these elements, as~$\vu=\du_1$, $\ru=\hu_2$, and~$\su=\hu_3$.

Since the coaction of~$G$ on~$\Au$ is defined in the same way as the coaction of~$G$ on~$A$, the elements~$\cu_4$ and~$\hu_j$ are group-like. Exactly as in Paragraph~\ref{CoalgBiprod}, we get
\begin{align*}
\Delta_{\Bu}(\cu_2) &= \frac{1}{2} \cu_2(\hu_1 + \hu_3) \ot \cu_2
+ \frac{1}{2} \cu_2 (\hu_1 - \hu_3) \ot \cu_3
\end{align*}
for the coproduct of~$\cu_2$ and
\begin{align*}
\Delta_{\Bu}(\du_1) &= \frac{1}{4} \du_1(\hu_1 + \hu_2 + \hu_3 + \hu_4) \ot \du_1
+ \frac{1}{4} \du_1 (\hu_1 + \zeta^2 \hu_2 - \hu_3 - \zeta^2 \hu_4) \ot \du_2 \\
&\quad + \frac{1}{4} \du_1 (\hu_1 - \hu_2 + \hu_3 - \hu_4) \ot \du_3
+ \frac{1}{4} \du_1 (\hu_1 - \zeta^2 \hu_2 - \hu_3 + \zeta^2 \hu_4) \ot \du_4 
\end{align*}
for the coproduct of~$\du_1$. As in Paragraph~\ref{CoalgBiprod}, the coproduct is also determined by its values on the generators:
\begin{lem}
We have $\Delta_{\Bu}(\ru) = \ru \ot \ru$ and $\Delta_{\Bu}(\su) = \su \ot \su$. Furthermore, we have
\begin{align*}
\Delta_{\Bu}(\uu) &= \frac{\K}{2} (\uu \ot \uu + \uu \ot \uu^3 +  \uu^3 \su \ot \uu - \uu^3 \su \ot \uu^3)
\end{align*}
and
\begin{align*}
\Delta_{\Bu}(\vu) &= \frac{1}{4} \vu (1 + \ru + \su + \ru \su) \ot \vu
+ \frac{1}{4} \vu (1 - \zeta^2 \ru - \su + \zeta^2 \ru \su) \ot \uu \vu \\
&\quad + \frac{1}{4} \vu (1 - \ru + \su - \ru \su) \ot \uu^2 \vu
+ \frac{1}{4} \vu (1 + \zeta^2 \ru - \su - \zeta^2 \ru \su) \ot \uu^3 \vu
\end{align*}
The counit is given on generators by~$\varepsilon_{\Bu}(\uu) = \varepsilon_{\Bu}(\vu) = \varepsilon_{\Bu}(\ru) = \varepsilon_{\Bu}(\su) = 1$.
\end{lem}

\subsection[The isomorphism between the first and the second case]{} \label{Isom}
The Hopf algebra~$\Bu$ is a Radford biproduct by construction, and as in the case of~$B$, it is reasonable to ask whether, and in which ways, it can also be written as an extension. Fortunately, it is not necessary to repeat this analysis for~$\Bu$, because~$\Bu$ is isomorphic to~$B$, as we show now. We begin with the following lemma:
\begin{lem}
There is an algebra homomorphism $f_{A'} \colon \Au \to B$ with the property that
\[f_{A'}(\xu) = u \qquad \text{and} \qquad f_{A'}(\yu) = v r\]
\end{lem}
\begin{pf}
We have to check the defining relations of~$\Au$, using Proposition~\ref{AlgBiprod}. It is obvious that $f_{A'}(\xu)^4 = 1$. For the second relation, we have
\begin{align*}
f_{A'}(\yu) f_{A'}(\xu)^3 &= v r u^3 = v u r = u v r = f_{A'}(\xu) f_{A'}(\yu)
\end{align*}
For the third relation, we have
\begin{align*}
f_{A'}(\yu)^2 &= v r v r = v (u^3 v r) r =  v^2 u^3
= \frac{\K}{2}(\ub + \zeta u + u^2 - \zeta u^3) u^3 \\
&= \frac{\K}{2}(u^3 + \zeta \ub + u - \zeta u^2)
= \frac{\K}{2}(\zeta \ub + f_{A'}(\xu) - \zeta f_{A'}(\xu)^2 + f_{A'}(\xu)^3)
\end{align*}
as required.
\qed
\end{pf}

The algebra homomorphism~$f_{A'}$ can be extended to the entire biproduct:
\begin{prop}
The map~$f \colon \Bu = \Au \star H \to B= A \star H$ given by
\[f(\au \star h) = f_{A'}(\au) (\ua \star h)\]
is an algebra isomorphism.
\end{prop}
\begin{pf}
We apply the universal property of a smash product stated in Corollary~\ref{UnivPropCross}. For the map~$f_H$ required there, we use the canonical embedding
\[f_H \colon H \to A \star H,~h \mapsto \ua \star h \]
of~$A$ into the biproduct~$B$. We have to check that
\[f_{A'}(h_\1.\au) f_H(h_\2) = f_H(h) f_{A'}(\au)\]
As pointed out after that corollary, we only need to check this on the algebra generators~$g_2$ and~$g_3$ of~$H$. In both cases, the condition can be written in the form
\begin{align*}
f_{A'}(g_i.\au) = f_H(g_i) f_{A'}(\au) f_H(g_i)^{-1}
\end{align*}
In this form, both sides depend multiplicatively on~$\au$, so that it is sufficient to verify this equation for $\au = \xu$ and~$\au = \yu$. Hence we have to consider four cases:
\begin{parlist}
\item
The case $h=g_2$, $\au=\xu$: We then have
\begin{align*}
f_H(g_2) f_{A'}(\xu) f_H(g_2)^{-1} &= r u r^{-1} = u^3 = f_{A'}(\xu^3) = f_{A'}(g_2.\xu)
\end{align*}

\item
The case $h=g_2$, $\au=\yu$: We then have
\begin{align*}
f_H(g_2) f_{A'}(\yu) f_H(g_2)^{-1} &= r (vr) r^{-1} = rv = u^3 v r = f_{A'}(\xu^3 \yu) = f_{A'}(g_2.\yu)
\end{align*}

\item
The case $h=g_3$, $\au=\xu$: We then have
\begin{align*}
f_H(g_3) f_{A'}(\xu) f_H(g_3)^{-1} &= s u s^{-1} = u = f_{A'}(\xu) = f_{A'}(g_3.\xu)
\end{align*}

\item
The case $h=g_3$, $\au=\yu$: We then have
\begin{align*}
f_H(g_3) f_{A'}(\yu) f_H(g_3)^{-1} &= s v r s^{-1} = u^2 v s r s^{-1} = u^2 v r
= f_{A'}(\xu^2 \yu) = f_{A'}(g_3.\yu)
\end{align*}
\end{parlist}
Because the elements~$u$, $vr$, $r$, and~$s$ are contained in its image, $f$~is surjective, and therefore bijective by dimension considerations.
\qed
\end{pf}

We note that, on the generators of~$\Bu$, the map~$f$ takes the values
\[f(\uu) = u \qquad \qquad f(\vu) = vr \qquad \qquad f(\ru) = r \qquad \qquad f(\su) = s \]
These equations can be used to give an alternative proof of the result above, which is based on Proposition~\ref{YetDrinfSec}. They are also used in the proof of the following corollary:
\begin{cor}
The map~$f$ is a Hopf algebra isomorphism between~$B$ and~$\Bu$.
\end{cor}
\begin{pf}
We still have to check that~$f$ is a coalgebra homomorphism, which amounts to the
conditions~$\db \circ f = (f \ot f) \circ \Delta_{\Bu}$ and~$\eb \circ f = \varepsilon_{\Bu}$. It is sufficient to check both conditions on the generators, where the second condition is obvious. The first condition is also obvious on the generators~$\uu$, $\ru$, and~$\su$, because their coproducts and the coproducts of the respective images are given essentially by the same formulas. For~$\vu$, we see by comparing Proposition~\ref{CoalgBiprod} and Lemma~\ref{YetDrinfSec} that
\[(f \ot f)(\Delta_{\Bu}(\vu)) = \db(v) (r \ot r) = \db(v) \db(r) = \db(vr) = \db(f(\vu))\]
which completes our proof, because a bialgebra homomorphism commutes automatically with the antipode.
\qed
\end{pf}


\end{document}